\documentclass[12pt]{amsart}
\usepackage{amsfonts}
\usepackage{amssymb}
\usepackage{amsmath}
\input xy
\xyoption {all}

\usepackage{amsthm}
\usepackage{amssymb}
\usepackage{latexsym}
\usepackage[dvips]{graphics}
\newtheorem{thm}{Theorem}[section] 
\newtheorem{cor}[thm]{Corollary}
\newtheorem{defn}[thm]{Definition}
\newtheorem{example}[thm]{Example}

\newtheorem{lemma}[thm]{Lemma}
\newtheorem{prop}[thm]{Proposition}
\newtheorem{remark}[thm]{Remark}

\DeclareMathOperator{\supp}{Supp}

\numberwithin{equation}{section}

\newcommand{\LL}{{\mathcal L}}
\newcommand{\U}{{\mathcal U}}
\newcommand{\A}{{\mathcal A}}

\newcommand{\CC}{{\mathcal C}}
\newcommand{\OO}{{\mathcal O}}

\newcommand{\Z}{\mathbb{Z}}
\newcommand{\N}{\mathbb{N}}
\newcommand{\Q}{\mathbb{Q}}

\newcommand{\C}{\mathbb{C}}

\newcommand{\PP}{\mathbb{P}}

\newcommand{\T}{\mathbb{T}}

\newcommand{\e}{\mathbf{e}}

\begin{document}
\date{}
\title{Multivariable Alexander invariants of hypersurface complements}

\author[Alexandru Dimca]{Alexandru Dimca}
\address{A.Dimca:  Laboratoire J.A. Dieudonn\'e, UMR du CNRS 6621,
                 Universit\'e de Nice-Sophia-Antipolis,
                 Parc Valrose,
                 06108 Nice Cedex 02,
                 FRANCE.}
\email
{dimca@math.unice.fr}

\author[Laurentiu Maxim ]{Laurentiu Maxim}
\address{L. Maxim : Department of Mathematics,
          University of Illinois at Chicago,
          851 S.Morgan, Chicago, Illinois, 60607
                 USA.}
\email
{lmaxim@math.uic.edu}

\subjclass[2000]{Primary
32S20, 32S22, 32S35, 32S60; Secondary
 14J70, 14F17, 14F45.
}

\keywords{hypersurface complement, Alexander polynomials, local system, Milnor fiber, logarithmic forms}

\begin{abstract}
We start with a discussion on Alexander invariants, and
then prove some general results concerning the divisibility of the
Alexander polynomials and the supports of the Alexander modules,
via Artin's vanishing theorem for perverse sheaves.
We conclude with explicit computations of twisted cohomology
following an idea already exploited in the hyperplane arrangement
case, which combines the degeneration of the Hodge to de Rham
spectral sequence to the purity of some cohomology groups.

\end{abstract}

\maketitle

\section{Introduction}
Alexander invariants in the form of Alexander modules, characteristic varieties
and Alexander polynomials have been recently intensively studied, in particular
in relation to the twisted cohomology of hypersurface arrangement complements, see for instance \cite{A},   \cite{COG}, \cite{CO},  \cite{CS}, \cite{Di2}, \cite{Di3},  \cite{Di4},
 \cite{DiLi},  \cite{ESV},  \cite{Li2},
 \cite{Li3}, \cite{Li8}, \cite{Max},  \cite{STV}.

In section \S2, after giving the basic definitions introducing the
{\it Alexander modules} $A^q(\U)$ of an affine hypersurface
arrangement complement $\U$, we investigate in
 Proposition \ref{prop2} the relation between {\it the first  non-trivial} Alexander polynomial in one variable and the corresponding Alexander polynomial in several
 variables. Proposition \ref{prop3} expresses the relation between the characteristic varieties defined using the Fitting ideals and  the characteristic varieties defined using the jumping loci of the cohomology with rank one local coefficients. Example \ref{ex1} treats the simplest {\it local} situations: the normal crossing case and the case of isolated non-normal crossing singularities, whose study was initiated by A. Libgober in \cite{Li8}.

In section \S3, Theorem \ref{thm1} relates the Alexander
invariants of the affine hypersurface arrangement complement
$\U=\C^{n+1} \setminus X$ to the  Alexander invariants of the
 complement $\U_{\infty}$ of the corresponding {\it link at infinity.} Theorems
\ref{thm2}, \ref{ucam} and Corollary \ref{cor1} estimate the support of the Alexander modules $A^q(\U)$ in terms of local properties of the projective closure $V= {\overline X   }$.

In section \S4, we recall and slightly extend the 
idea of combining the degeneration of the Hodge to de Rham spectral
sequence to the purity of some cohomology groups (used first by  Esnault, Schechtman and  Viehweg in
\cite{ESV} and  by Schechtman, Terao and Varchenko in \cite{STV}), see Corollary  \ref{Delcor} and
Proposition \ref{proppure}. Examples  \ref{linesandconic1} and
\ref{linesandconic2} illustrate this approach by looking at some
arrangements of lines and conics in the plane. Though these
examples may be treated using the results by Cogolludo in \cite{COG}, we feel
that our approach is more general and hence more likely to extend
to other situations.

In the last section we consider the complement $\U_0$ of an
arbitrary {\it projective} hypersurface arrangement $V$, and,
after a short general discussion, we revisit from a new
perspective a useful result by Randell saying what happens to the
twisted cohomology of a plane curve complement when we add an
extra line, see Corollary  \ref{curvearr}. Coming back to
dimension $n \geq 2$, Example \ref{onecomp} discusses the already
interesting case when $V$ is irreducible and has only isolated
singularities. This case leads, in particular, to examples where
for  $m=n, n+1$ and some  rank one local coefficients
$\LL_{\beta}$ on $\U_0$ one has
$$\dim  H^m(\U_0,\LL_{\beta}) >  \dim  H^m(\U_0,\C).$$
By the minimality property of hyperplane arrangement complements,
it is known that the above inequality is impossible for such type
of complements, \cite{DP}. We conclude by a detailed study of the
case when $V$ has two  irreducible components, each of them having
only isolated singularities.

 Throughout the paper we usually work
with complex coefficients $\C$, although the study of finite field coefficients
is very important, due for instance to torsion open questions, see 
\cite{CDS}, \cite{MS}. Our choice is imposed by the analytic tools used in the last two sections. Most of the results in the previous sections hold over arbitrary fields.

\section{Multivariable Alexander invariants}
\subsection{Algebraic Preliminaries}
Let $R$ be a commutative ring with unit, which is  Noetherian and
a unique factorization domain (e.g. the ring  of complex Laurent
polynomials in $s$ variables, $s \geq 1$). Let $A$ be a finitely
generated $R$-module, and $M$ a $(n \times m)$ presentation matrix
of $A$ associated to an exact sequence
$$ R^m \to R^n \to A \to 0 .$$
The $i$-th \emph{elementary ideal} $\mathcal{E}_i(A)$ of $A$ is
the ideal in $R$ generated by the $(n-i)\times (n-i)$ minor
determinants of $M$, with the convention that $\mathcal{E}_i(A)=R$
if $i \geq n$, and $\mathcal{E}_i(A)=0$ if $n-i>m$. Let
$\Delta_i(A)$ be the generator of the smallest principal ideal in
$R$ containing $\mathcal{E}_i(A)$, i.e. the greatest common
divisor of all elements of $\mathcal{E}_i(A)$. $\Delta_i(A)$ is
called the $i$-th \emph{characteristic polynomial} of $A$. Note
that $\Delta_{i+1}(A)$ divides $\Delta_i(A)$ in $R$ for all $i$
since $\mathcal{E}_i(A) \subset \mathcal{E}_{i+1}(A)$. In
particular, if $R$ is a principal ideal domain (e.g. the ring  of
complex Laurent polynomials in one variable), then
$\mathcal{E}_i(A)$ is a principal ideal  generated exactly  by
$\Delta_i(A)$.

As an example, for any ring $R$, assume that $A=R^s \oplus R/(\lambda_1) \oplus \cdots \oplus
R/(\lambda_r)$, where $\lambda_j$ ($j=1,2,\cdots,r$) are non-zero
elements in $R$ such that $\lambda_{j+1}|\lambda_j$. Then we have
$\Delta_i(A)$ is $0$, $\lambda_{i-s+1} \cdots \lambda_r$, or $1$,
according to whether $0 \leq i \leq s-1$, $s\leq i \leq s+r-1$, or
$s+r \leq i$.

The \emph{support} $\text{Supp}(A)$ of $A$ is the reduced
sub-scheme of $\text{Spec}(R)$ defined by (the \emph{order ideal})
$\mathcal{E}_0(A)$. Since
$$\sqrt { \mathcal{E}_0(A)}= \sqrt {Ann(A)}$$
this is the usual notion of support in algebraic geometry based on
the \emph{annihilator ideal} $\text{Ann}(A)$ of the module $A$. In
particular, for a prime ideal $P \subset R$, $P \in
\text{Supp}(A)$ if and only if the localized module $A_P$ is
non-zero.

The support $\text{Supp}(A)$ is also called the
\emph{first characteristic variety} of $A$, and we define the
$i$-th \emph{characteristic variety} $V_i(A)$ of $A$ to be the
reduced sub-scheme of $\text{Spec}(R)$ defined by the ($i$-th
\emph{Fitting ideal}) ideal $\mathcal{E}_{i-1}(A)$.

\noindent Note that $\text{codim}V_i(A)>1$ implies $\Delta_{i-1}
(A)=1$, i.e. the corresponding Alexander polynomial carries no
information.

All definitions above are independent (up to multiplication by a
unit of $R$) of the choices involved, thus the characteristic
varieties and polynomials of $A$ are invariants of the
$R$-isomorphism type of $A$.

We state for future reference the following 'divisibility'
properties of the polynomials and characteristic varieties (for
proofs, see \cite{S} and \cite{Li7}):
\begin{lemma}\label{div}
\ \\

\begin{itemize}
\item If $A$, $B$ are finitely generated $R$-modules, then:
$\Delta_0(A \oplus B)=~\Delta_0(A) \times \Delta_0(B)$. \item If
$A$ and $B$ are finitely generated $R$-modules then:
\begin{center}
$  \supp (A \otimes_R B)= \supp   (A) \cap \supp   (B)$.
\end{center}
\item If $A$ is a submodule of $B$, then for all $i$,
$\Delta_i(A)$ divides $\Delta_i(B)$. \item If $0 \to A \to B \to C
\to 0$ is a short exact sequence of finitely generated
$R$-modules, then the following hold:
\begin{enumerate}
\item $\Delta_0(B)=\Delta_0(A)\times \Delta_0(C)$; \item for all
$i$, $\Delta_i(B)$ divides $\Delta_i(A)\times \Delta_0(C)$; \item
If $\Delta_0(C)=1$, then $\Delta_i(A)=\Delta_i(B)$ for all $i$;
\item $ \supp  (B)= \supp  (A) \cup  \supp  (C)$; \item
For $i \geq 2$: $V_i(C) \subset V_i(B) \subset V_i(C) \cup
(V_{i-1}(C) \cap  \supp (A))$.
\end{enumerate}
\end{itemize}
\end{lemma}

\bigskip

\subsection{Alexander Invariants of Hypersurface Complements}\label{def}
Let $V$ be a reduced hypersurface in $\mathbb{CP}^{n+1}$, defined
by a homogeneous equation: $f=f_1\cdots f_s =0$, where $f_i$ are
the irreducible factors of $f$, and $V_i=\{f_i=0\}$ the
irreducible components of $V$. We fix a hyperplane $H$ in
$\mathbb{CP}^{n+1}$ which we call 'the hyperplane at infinity'.
Let $\mathcal{U}$ be the (affine) hypersurface complement
$\mathcal{U}=\mathbb{CP}^{n+1} \setminus (V \cup H)$.
(Alternatively, $\U$ may be regarded as the complement of a
hypersurface in the affine space $\C^{n+1}$.) Then $H_1
(\mathcal{U}) \cong \mathbb{Z}^s$ (\cite{Di}, (4.1.3), (4.1.4)),
generated by the meridian loops $\gamma_i$ about the non-singular
part of each irreducible component $V_i$, for $i=1,\cdots, s$. If
$\gamma_\infty$ denotes the meridian about the hyperplane at
infinity, then in $H_1(\mathcal{U})$ there is a relation:
$\gamma_\infty + \sum {d_i \gamma_i} = 0$, where $d_i=deg(V_i)$.

Note that $\mathcal{U}$ is affine, therefore has the homotopy type
of a finite CW complex. Let $\mathcal{U}^{ab}$ be the universal
abelian cover of $\mathcal{U}$, i.e. the covering associated to
the commutator subgroup of $\pi_1(\mathcal{U})$, or equivalently,
the covering associated to the kernel of the linking number
homomorphism $lk : \pi_1(\mathcal{U}) \to \mathbb{Z}^s$, which
maps a loop $\alpha$ to $(\text{lk} (\alpha, V_1 \cup
-d_1H),\cdots, \text{lk} (\alpha, V_s \cup -d_sH))$. The group of
covering transformations of $\mathcal{U}^{ab}$ is isomorphic to
$\mathbb{Z}^s$ and acts on the covering space. By choosing fixed
lifts of the cells of $\mathcal{U}$ to $\mathcal{U}^{ab}$, we
obtain a free basis for $C_\ast$, the cellular cell complex of
$\mathcal{U}^{ab}$, as a $\mathbb{Z}[\mathbb{Z}^s]$-module. The
isomorphism determined by the meridians $\{\gamma_i\}$ enables us
to identify $\mathbb{Z}[\mathbb{Z}^s]$ with
$\mathbb{Z}[t_1,t_1^{-1},\cdots, t_s,t_s^{-1}]$, the ring of
integral Laurent polynomials in $s$ variables. When $s=1$ we set
$t_1=t$.

For reasons that will become transparent later, our base ring will
always be the ring of complex Laurent polynomials in $s$
variables, $\mathbb{C}[t_1,t_1^{-1},\cdots, t_s,t_s^{-1}]$, which
we denote by $R_s$. Note that $R_s$ is a regular Noetherian
domain, and in particular it is factorial. As a groups ring, $R_s$
has a natural involution denoted by an overbar, sending each $t_i$
to $\bar t_i :=t_i ^{-1}$. To an $R_s$-module $A$, we associate
the conjugate $R_s$-module, still denoted by $A$, with the same
underlying abelian group but with the $R_s$-action given by $(r,a)
\mapsto {\bar r}\cdot a$, for $a \in A$ and $r \in R_s$.

\begin{remark} \label{rem0} \rm

Though the ring  $R_s$ is commutative, it should be regarded as a
quotient ring of $\C [\pi_1(\mathcal{U})]$, which is
non-commutative in general. Because of that, one should be careful
to distinguish the right from the left $R_s$-modules. If, for
instance, $A$ is a left $R_s$-module, then the associated right
$R_s$-module is the module conjugate to $A$, whose module
structure is given by
$$a \cdot r:= {\bar r} \cdot a$$
for all $a \in A$ and $r \in R_s$. This corresponds to regarding
any left $\C [\pi_1(\mathcal{U})]$-module $A$ as a right $\C
[\pi_1(\mathcal{U})]$-module by setting $a \cdot \gamma = \gamma
^{-1} \cdot a$, for all $a \in A$ and $\gamma \in
\pi_1(\mathcal{U})$, and extending by linearity. Following
\cite{DK}, p. 97, we regard in this paper $C_\ast ^0= C_\ast
\otimes \mathbb{C}$ as a complex of right $R_s$-modules.

\end{remark}

Define a local coefficient system $\mathcal{L}$ on $\mathcal{U}$,
with stalk $R_s$ and action of a loop $\alpha \in
\pi_1(\mathcal{U})$ determined by (left) multiplication by
$\prod_{j=1}^s {(t_j)^{\text{lk}(\alpha,V_j \cup -d_jH)}}$. In
particular, the action of the meridian $\gamma_i$ is given by
multiplication by $t_i$. Let $\mathcal{L}^{\vee}$ be the dual
local system, whose stalk at a point $y \in \U$ is
$\mathcal{L}^{\vee} _y:= \text{Hom}(\mathcal{L}_y,R_s)$, and let
$\alpha \in \pi_1(\mathcal{U},y)$ act on $\varphi \in
\mathcal{L}^{\vee} _y$ by:
$$(\alpha \cdot \varphi) (m) := \varphi (\alpha ^{-1} \cdot m) \ , \ \ m \in
\mathcal{L}_y.$$ We denote by $\bar {\mathcal{L}}$ the local
system obtained from $\mathcal{L}$ by composing all module
structures with the involution of $R_s$ (i.e. by changing the
stalks of $\mathcal{L}$ from left into right $R_s$-modules). The
perfect pairing $$\bar {\mathcal{L}} \otimes_{R_s} \mathcal{L} \to
R_s$$ given by
$$(f,g) \mapsto {\bar f} \cdot g$$ on the stalk over a basepoint,
tell us that there is an isomorphism of local systems on
$\mathcal{U}$: $$\mathcal{L}^{\vee} \simeq \bar {\mathcal{L}}.$$

The \emph{ universal homology k-th Alexander invariant } $A_k(\U)$
of $\U$ is by definition the $R_s$-module $H_k(C_\ast ^0)$, or
equivalently $H_k(\mathcal{U};\mathcal{L})$. This is the group
$H_k(\mathcal{U}^{ab};\mathbb{C})$ considered as a $R_s$-module
via the covering transformations (see \cite{Hat}, Example 3H.2).
Similarly, the \emph{ universal cohomology k-th Alexander
invariant } $A^k(\U)$ of $\U$ is by definition the $k$-th
cohomology module of the dual complex $\text{Hom}_{R_s}(C_\ast ^0,
{R_s})$. Here $R_s$ is considered with the induced right
$R_s$-module structure as explained in Remark \ref{rem0}. Based on
our previous considerations on local systems, $A^k(\U)$ is just
$H^k(\mathcal{U};\mathcal{L}^{\vee})$. This may be also regarded
as the $k$-th cohomology with compact support and complex
coefficients of $\mathcal{U}_b ^{ab}$, where $\U_b$ is the compact
manifold with boundary obtained from $\mathbb{CP}^{n+1}$ by
removing a small open regular neighborhood of the divisor $V \cup
H$ (compare \cite{Hat}, Prop. 3H.5).

Note that, since $\mathcal{U}$ is a $(n+1)$-dimensional affine
variety, the modules $A^k(\U)$ and resp. $A_k(\U)$ are trivial for
$k > n+1$. Moreover, since the stalks of $\mathcal{L}$ are
torsion-free, $A_{n+1}(\U)$ is also a torsion-free $R_s$-module
(see \cite{Sh}, Example 6.0.6).

As in the classical knot theory, by using a deformation retract
argument, one could  define the universal abelian invariants above
after replacing $\U$ by the manifold with boundary  $\U_b$,
obtained from $\mathbb{CP}^{n+1}$ by removing a small open regular
neighborhood of the divisor $V \cup H$. Now, since the chain
complex $C_\ast (\U_b ^{ab})$ is of finite type, and since $R_s$
is Noetherian, this implies that all these universal Alexander
modules are finitely generated. Hence their characteristic
varieties and polynomials are well-defined. The associated
characteristic varieties, in particular the supports, become
sub-varieties of the $s$-dimensional torus
$\T^s=(\mathbb{C}^\ast)^s$, which is regarded as the set of closed
points in $\text{Spec}(R_s)$. More precisely, for  $\lambda
=(\lambda _1,\cdots, \lambda  _s) \in \T^s$, we denote by
$m_{\lambda}$ the corresponding maximal ideal in $R_s$ and by
$\C_{\lambda}$ the quotient $R_s/m_{\lambda}R_s$. This quotient is
isomorphic to $\C$ and the canonical projection
\begin{equation} \label{eval}
 \rho _{\lambda}:R_s \to R_s/m_{\lambda}R_s=\C_{\lambda}
\end{equation}
corresponds to replacing $t_j$ by $\lambda   _j$ for $j=1,...,s.$
Here we regard $\C_{\lambda}$ as a (left) $R_s$-module, with an
involution given by the complex conjugation (which is compatible
with the one induced from $R_s$ since $\lambda_j \in \T^1$).

If $A$ is an $R_s$-module, we denote be $A_{\lambda}$ the
localization of $A$ at the maximal ideal $m_{\lambda}$. For
$A=R_s$, we use the simpler notation $R_{\lambda}$ when there is
no danger of confusion.   If $A$ is of finite type, then $A=0$ if
and only if $A_{\lambda}=0$ for all   $\lambda  \in \T^s$. More
precisely
$$ \text{Supp}(A)= \{ \lambda  \in \T^s; A_{\lambda} \ne 0\}$$
In particular $A_0(\U)=\C_{\bf 1}$, where ${\bf 1}=(1,\cdots, 1)$ and hence
\begin{equation} \label{supp1}
 \text{Supp}(A_0(\U))  = \{{\bf 1}\}.
\end{equation}
We denote by $V_{i,k}(\mathcal{U})$ the $i$-th characteristic
variety associated to the homological Alexander module
$A_k(\U)$, and similarly denote by
$\Delta_{i,k}(\mathcal{U})$ the associated characteristic
polynomials. The notations $V^{i,k}(\mathcal{U})$ and $\Delta^{i,k}(\mathcal{U})$
denote the similar objects associated to the cohomological Alexander invariants
$A^k(\U)$.

\subsection{Homology versus Cohomology Alexander Modules}

It is natural to ask what are the relations between the homology and the cohomology universal Alexander modules. Or to find the relations between $V_{i,k}(\mathcal{U})$ and  $V^{i,k}(\mathcal{U})$;
and between $\Delta_{i,k}(\mathcal{U})$ and  $\Delta^{i,k}(\mathcal{U})$.

\bigskip

Some answers to this question can be  given as follows. The
cohomology modules may be related to the homology modules by the
Universal Coefficient spectral sequence (see \cite{Hi}, p.20 or
\cite{L}, Thm. 2.3).
\begin{equation} \label{spsq1}
\text{Ext}^q_ {R_ {s}}  (A_p(\mathcal{U}),R_s) \Rightarrow
A^{p+q}(\mathcal{U}).
\end{equation}
Using the exactness of the localization (see \cite{W}, p. 76), we
get the following spectral sequence for any $ \lambda  \in \T^s$.
\begin{equation} \label{spsq2}
\text{Ext}^q_{R_ {\lambda}} (A_p(\mathcal{U})_{\lambda} ,R_
{\lambda} ) \Rightarrow A^{p+q}(\mathcal{U})_{\lambda}.
\end{equation}
For a fixed $ \lambda  \in \T^s$, we define
\begin{equation} \label{defk}
k( \lambda) =\text{min} \{m \in \N; A_m(\mathcal{U})_{\lambda} \ne
0\}.
\end{equation}
Then the spectral sequence \ref{spsq2} implies the following.

\begin{prop} \label{prop1}

For any $ \lambda  \in \T^s$, $A^{k}(\U)_ {\lambda}=0$ for $k< k( \lambda)$ and
\begin{equation} \label{homcohom}
A^{k( \lambda)}(\U)_ {\lambda}  =\text{Hom}(A_{k( \lambda)}(\U)_
{\lambda},R_ {\lambda})
\end{equation}

\end{prop}

This equality shows in particular that one may have $A^{k( \lambda)}(\U)_ {\lambda}=0$,
even when $A_{k( \lambda)}(\U)_ {\lambda} \ne 0$, e.g. when the last module is torsion,
which is often the case, e.g. see \ref{supp1}.

\subsection{Multivariable versus one variable Alexander Modules}

Consider a family of integral weights $\e=(e_1,\cdots,e_s) \in \mathbb{Z}^s$, and let
$$q:=\text{g.c.d.}(e_1,\cdots,e_s).$$
Consider the morphism $p(\e):R_s \to R_1$ defined by $t_i \mapsto
t^{e_i}$, inducing a (left) $R_s$-module structure on $R_1$. Let
$\LL (\e)$ be the local system on $\U$ with stalk $R_1$ and
monodromy action for a loop $\alpha \in \pi_1(\mathcal{U})$ given
by multiplication by $t^{\sum e_j\text{lk}(\alpha,V_j \cup
-d_jH)}$.

The corresponding homology groups $H_k(\U,\LL (\e))=H_k(C_\ast^0 \otimes _{R_s}R_1)$ are
finite type $R_1$-modules, and hence they have associated characteristic
varieties $V_{i,k}(\mathcal{U},\e)$ and
Alexander polynomials
$\Delta_{i,k}(\mathcal{U},\e  )$.

\bigskip

It is natural to ask under which conditions the equalities
$$ \Delta_{i,k}(\mathcal{U},\e  )(t)=(t^q-1)\Delta_{i,k}(\mathcal{U})(t^{e_1},\cdots,t^{e_s})$$
do hold? Something like this works in classical knot theory, more
precisely for oriented multilinks in $S^3$ with at least 2
components, where the case $i=0$, $k=1$ is considered (see
\cite{EN}, Prop. 5.1, and also \cite{Mi}, Lemma 10.1 for the case
of weight $(1,\cdots,1)$).

For the  weight ${\bf 1}=(1,1,...,1)$, we call the coresponding Alexander polynomials the
usual (or, univariable ) Alexander polynomials and we denote them by $\Delta_{i,k}^T(\mathcal{U})$
(see below for some explanation).

If the equality in  Question 2 holds for all but finitely many
multi-indices $\e$, then the 1-variable polynomials $
\Delta_{i,k}(\mathcal{U},\e  )$ determine (up-to a unit in $R_s$)
the multi-variable polynomial $\Delta_{i,k}(\mathcal{U})$ (see
\cite{Ci}, Lemma 2.2).

\bigskip

Some insight into this question can be obtained as follows. We
consider only  the simplest case, namely $\e= {\bf 1}$, and leave
the other cases to the interested reader.

Note that the universal abelian covering $\U^{ab} \to\U$ corresponds to the kernel $ K^{ab}$ of the abelianization morphism
$$ \pi _1(\U) \to H _1(\U).$$
The total linking number covering $\U^{T} \to\U$ corresponds to the kernel $ K^T$ of the morphism
$$ \pi _1(\U) \to H _1(\U)=\Z^s \to \Z$$
where the second morphism is $\sum c_j \gamma_j \mapsto \sum c_j$.
It follows that $\U^{ab} \to\U^T$ is a covering with deck transformation group $G=K^T/ K^{ab}$
identified to the subgroup
$$\{c \in \Z^s; \sum c_j=0\}.$$
The complex $C_\ast^0$ is a complex of free $R_s$-modules of finite rank and the derivatives are  $R_s$-liniar.
It follows that we can regard this complex as being a complex $\CC_\ast^0$ of free $\OO_{\T^s}$-modules
on the affine variety $\T^s$.

Since $\U^T=\U^{ab}/G$, it follows that the complex of singular chains of  $\U^T$ is
\begin{equation} \label{homcov}
C_*(\U^T)=C_*(\U^ {ab}  )_G=(C_\ast^0)_G
\end{equation}
(see \cite{W}, p.204). Here
\begin{equation} \label{definv}
(C_p^0)_G=C_p^0/<gm-m;~ g \in G, ~ m \in C_p^0>.
\end{equation}
Using the fact that the group $G$ is generated by the elements
having an 1 as the $i$-th coordinate, a $-1$ as the $j$-th
coordinate (for $i<j$) and all the other coordinates zero, we see
that $(C_p^0)_G$ is the quotient of $C_p^0$ by the submodule
$$<(t_i-t_j)m;~ m \in C_p^0>.$$
It follows that the associated sheaf $(\CC_p^0)_G$ is just the
restriction (as a coherent sheaf) of $\CC_p^0$ to the
1-dimensional subtorus $S=\{(t,t,...,t) \in \T^s \}$, i.e.
$(\CC_p^0)_G= \CC_p^0 \otimes _{\OO_{\T^s }}\OO_{S }.$
Unfortunately, the inclusion $S \to \T^s$ is not a flat morphism
(see \cite{Ha}, p. 254), and hence the restriction to $S$ does not
commute to taking homology.\newline However, by our discussion
above
$$(C_p^0)_G=C_p^0 \otimes _{R_s}R_1,$$ with the (left) $R_s$-module
structure on $R_1$ induced by $p(\bf 1)$. Use now the K\"unneth
spectral sequence (see \cite{W}, p.143), and get
\begin{equation} \label{spsq3}
E^2_{p,q}=Tor_p^{R_s}(A_q(\U),R_1)  \Rightarrow H_{p+q}((C_\ast^0)_G)=  A_{p+q}^T(\U).
\end{equation}
For $a \in \T^1 =S=\{(t,t,...,t) \in \T^s \}$, we get by localization a new  K\"unneth spectral sequence,
namely
\begin{equation} \label{spsq4}
E^2_{p,q}=Tor_p^{R_a}(A_q(\U)_a,R_{1,a})  \Rightarrow H_{p+q}((C_\ast^0)_G)_a.
\end{equation}
In particular we get the following.

\begin{prop} \label{prop2}

For any $ a  \in \T^1$, $A_{k}^T(\U)_ {a}=0$ for $k< k(a)$ and
\begin{equation} \label{loc2}
A_{k(a)}(\U)_a \otimes_{R_a}R_{1,a}=A_{k(a)}^T(\U)_a
\end{equation}
In particular, for any $ a \in \T^1 =S$, the multiplicity of the root $t=a$
in the polynomials $ \Delta_{i,k(a)}^T(\mathcal{U})(t)$ and $ \Delta_{i,k(a)}(\mathcal{U})(t, \cdots,t)$
is the same.

\end{prop}

\proof To get the second claim, note that any presentation
$$ R_a^m \to R_a^n \to A_{k(a)}(\U)_a \to 0$$
yields by tensor product a presentation
$$ R_{1,a}  ^m \to R_ {1,a} ^n \to A_{k(a)}^T(\U)_a \to 0.$$
\endproof

\subsection{Characteristic varieties as jumping loci of  rank-1 local systems}

Let $\lambda =(\lambda   _1,\cdots, \lambda  _s) \in \T^s$ and
denote by $\LL_{ \lambda} $ the local coefficient system on $\U$
with stalk $\C = \C_{\lambda}$ and  action of a loop $\alpha \in
\pi_1(\mathcal{U})$ determined by multiplication by $\prod_{j=1}^s
{( \lambda _j)^{\text{lk}(\alpha,V_j \cup -d_jH)}}$. We let
$\mathcal{L}^{\vee}_{\lambda} \simeq \LL_{\lambda ^{-1}}$ be the
dual local system, where $\lambda ^{-1}:= (\lambda   _1
^{-1},\cdots, \lambda _s ^{-1}) \in \T^s$.

One can define new  \emph{topological characteristic varieties}  by setting
$$V_{i,k}^t(\mathcal{U})=\{\lambda  \in \T^s; \text{dim}H_k(\U,\LL_{ \lambda}  )> i\}$$
and
$$V^{i,k}_t(\mathcal{U})=\{\lambda  \in \T^s; \text{dim}H^k(\U,\LL_{ \lambda}  )> i\}.$$

It is natural to investigate the relations between the two types of characteristic varieties. Some cases are considered in \cite{Li7}, \cite{Li8}.

\bigskip

Here is a general approach to this question. It is known that
$$ H_k(\U, \LL_{ \lambda})=H_k(C_{\ast}^0 \otimes _{R_s}\C_{ \lambda}).$$
Using the K\"unneth spectral sequence, we get
\begin{equation} \label{spsq5}
E^2_{p,q}=Tor_p^{R_s}(A_q(\U),\C_{ \lambda}  )  \Rightarrow H_{p+q}(\U, \LL_{ \lambda}).
\end{equation}
Now since the localization is exact, the base change for Tor under
$R_s \to R_{ \lambda}$ (see  \cite{W}, p. 144), yields a new
spectral sequence
\begin{equation} \label{spsq6}
E^2_{p,q}=Tor_p^{R_ { \lambda}  }(A_q(\U)_{ \lambda}  ,\C_{ \lambda}  )  \Rightarrow H_{p+q}(\U, \LL_{ \lambda}).
\end{equation}
This proves the first claim of the next result.

\begin{prop} \label{prop3} For any point $\lambda  \in \T^s$, one has the following.

\medskip

\noindent (i) $\min \{m\in \N, ~ H_{m}(\U, \LL_{
\lambda})\ne 0 \}=\min \{m\in \N,~ \lambda \in
 \supp   (A_m(\U))
 \}=k( \lambda) .$

\medskip

\noindent (ii) $\dim H_{k( \lambda)}(\U, \LL_{
\lambda})=\max \{m\in \N,~ \lambda \in V_{m,k( \lambda)}(\U)
\} $.

\end{prop}

\proof To prove the second claim, note that the  spectral sequence  \ref{spsq6} yields
$$H_{k( \lambda)}(\U, \LL_{ \lambda})=A_{k( \lambda)}(\U)_{ \lambda}/m_{ \lambda}A_{k( \lambda)}(\U)_{ \lambda}.$$
Let $n$ be the dimension of these two vector spaces. Then by
Nakayama's Lemma, the module $A_{k( \lambda)}(\U)_{ \lambda}$ is
generated by $n$ elements over the local ring $R_{ \lambda}$. In
other words, there is presentation
$$ R_ { \lambda}  ^m \to R_ { \lambda}   ^n \to A_{k({ \lambda}  )}(\U)_ { \lambda} \to 0.$$
Moreover, the first morphism is given by a matrix $M$ whose entries $m_{ij}$ are all in the maximal ideal
$m_{ \lambda}$. The second claim follows now by the definition of the characteristic varieties.

\begin{remark} \label{rem1} \rm
Note that there is also a spectral sequence
\begin{equation}\label{ss}
E_2^{p,q}=\text{Ext}^q_ {R_ { \lambda}  } (A_p(\mathcal{U})_{
\lambda}  ,\C_ { \lambda}  ) \Rightarrow H^{p+q}(\mathcal{U},
\LL_{ {\lambda}^{-1}}     ).
\end{equation}
Here $\C_ { \lambda}$ is considered with the right $R_s$-module
structure as indicated in remark \ref{rem0}. This is why in the
abutment of the spectral sequence \ref{ss}, we obtain cohomology
with coefficients in the dual local system $\LL ^{\vee} _{\lambda}
\simeq \LL_{ {\lambda}^{-1}}$.\newline The above spectral sequence
yields that $H^m(\mathcal{U}, \LL_{ {\lambda}^{-1}})=0$ for $m<k(
\lambda)$ and $H^{k( \lambda)} (\mathcal{U}, \LL_{
{\lambda}^{-1}})=Hom_{R_{ \lambda}}(A_ {k( \lambda)}
(\mathcal{U})_{ \lambda}, \C_{ \lambda})$. However
$$Hom_{R_{ \lambda}}(A_ {k( \lambda)}  (\mathcal{U})_{ \lambda}, \C_{ \lambda})=
Hom_{\C}(A_ {k( \lambda)}  (\mathcal{U})_{ \lambda}/m _{ \lambda}
A_ {k( \lambda)}  (\mathcal{U})_{ \lambda}           , \C_{
\lambda})$$ and hence
\begin{equation} \label{dualls}
H_{k( \lambda)}(\mathcal{U}, \LL_{ \lambda})^*=H^{k( \lambda)}(\mathcal{U}, \LL_{ {\lambda}^{-1}}),
\end{equation} 
(compare \cite{Di2}, p.50 and p. 69). The case $k=1$ of this useful formula was
established in \cite{MS}, Remark 5.2. Note that this formula holds over arbitrary fields,
with the same proof as above.

\end{remark}

\begin{remark} \label{remlocal} \rm
All the results in this section so far hold for the local setting
as well, i.e. when $\U$ is the complement of a hypersurface germ
in a small ball. The first part of the  example below corresponds to the germ of a
normal crossing divisor. The second part of the  example below corresponds to isolated
non-normal crossing divisors (for short INNC), see \cite{DiLi}, \cite{Li8},
\cite{Li9}.

Similarly, instead of localizing at a point, one may localize along the hyperplane at infinity $H$, i.e. replace  $\U$ by
 $\U_{\infty}=\U\cap S_{\infty}    $, where $ S_{\infty}$ is a large enough sphere in $ \C^{n+1}$, see Theorem \ref{thm1} below.
\end{remark}

\begin{example} \label{ex1}  \rm

\noindent (i) Let $\U=(\C^*)^s \times \C^{n+1-s}$ for some integer $0 \le s \le
n+1$. Then the universal abelian covering $\U^{ab}$ is
contractible and then $A_0(\U)=\C_{\bf 1}$ and  $A_k(\U)=0$ for
$k>0$. Therefore, by the spectral sequence \ref{spsq1} we get
$A^k(\U)\cong \text{Ext}^k_{R_s}(\C_{\bf 1}, R_s)$ for all $k \geq
0$. Using the free resolution of $\C_{\bf 1}$ given by the Koszul
complex of the regular sequence $\{x_j=t_j-1\}_{ j=1,...,s}$ in
the ring $R_s$ (\cite{W}, Cor. 4.5.5), we obtain that $A^k(\U) =0$
for $k \ne s$ and $A^s(\U)=\C_{\bf 1}$ (\cite{W}, Ex. 4.5.2 and
Cor 4.5.4). Therefore, for any $ \lambda \ne {\bf 1}$, Proposition
\ref{prop1} shows that the corresponding cohomology Alexander
modules satisfy $A^k(\U)_{ \lambda}  =0$ for any $k$. Moreover
$H_k(\mathcal{U}, \LL_{ \lambda})=H^k(\mathcal{U}, \LL_{
{\lambda}^{-1}})=0$ for any  $k$.

\medskip

\noindent (ii) Let $(Y,0)$ be an INNC singularity at the origin of $ \C^{n+1}$.
Set $\U(Y,0)=B \setminus Y$, where $B$ is a small open ball centered at the origin in
$ \C^{n+1}$. Assume that $n \geq 2$. Then the universal abelian cover  $\U(Y,0)^{ab}$ of
 $\U(Y,0)$ is $(n-1)$-connected, see Libgober \cite{Li8}. More precisely, it is a bouquet of $n$-spheres,
see  \cite{DiLi}, and hence $A_0(\U(Y,0))=\C_{\bf 1}$ and
$A_k(\U(Y,0))=0$ for $k \ne n$. As in (i) above, we get
$A^k(\U(Y,0))\cong \text{Ext}^k_{R_s}(\C_{\bf 1}, R_s)$ for all $k
<n$. For $ \lambda \ne {\bf 1}$ this yields
${A^k(\U(Y,0))}_{\lambda}=0$ for $k<n$, and therefore $H^k(\U(Y,0)
, \LL_{ {\lambda}})=0$ for any  $k<n$.

\end{example}

\section{Divisibility Results and Characteristic Varieties}

In this section we give an algebraic-geometrical interpretation
for the multi-variable Alexander invariants of the hypersurface
complement, similar in flavor to the one-variable case described
in \cite{Max}, but see also the reformulation of these results in
\cite{Di4}. We will use an approach based on the general theory of
perverse sheaves, close to the one presented in \cite{Di4} (see
also \cite{COD} and \cite{Di2}). Note that the supports and
characteristic polynomials $\Delta_0$ of the multi-variable
Alexander modules are the analogue of the set of roots of the
Alexander polynomials and respectively Alexander polynomials in
the one-variable case (cf. \cite{Max}, \cite{Di4}).

The first result is an extension of \cite{Li6}, Theorem 3.2, to
arbitrary hypersurface singularities . Let $S_{\infty}$ be a
sphere of sufficiently large radius in $\mathbb{C}^{n+1}=
\mathbb{CP}^{n+1} \setminus H$ (or equivalently, the boundary of a
sufficiently small tubular neighborhood of $H$ in
$\mathbb{CP}^{n+1}$). Let $V_{\infty}=S_{\infty} \cap V$ be the
link of $V$ at infinity, and $\mathcal{U}_{\infty}=S_{\infty}
\setminus V_{\infty}$ its complement.

\begin{thm}\label{thm1}
For all $i$, and all $k \leq n$: $V_{i,k}(\mathcal{U}) \subset
V_{i,k}(\mathcal{U}_{\infty})$, and $\Delta_{i,k}(\mathcal{U}) |
\Delta_{i,k}(\mathcal{U}_{\infty})$. Moreover, for $k<n$, these inclusions
and divisibility conditions are replaced by equalities.
\end{thm}

\begin{proof} The case $n=1$ is considered in \cite{Li6}. In fact in this situation
one sets, for $i \leq 1$ and $k \leq 1$, $V_{i,k}(\mathcal{U}_{\infty})$ to be the $k$-th
characteristic variety of the $i$-th homology module of the covering space of $\mathcal{U}_{\infty}$
corresponding to the kernel of the composition
$$\pi_1(\mathcal{U}_{\infty}) \to \pi_1(\mathcal{U})\to H_1(\mathcal{U}).$$
For $n
\geq 2$, the theorem is an easy consequence of the Lefschetz
hyperplane theorem. Indeed, as in the proof of Theorem 4.5 of
\cite{Li2}, it follows that $\pi_1(\mathcal{U}) \cong
\pi_1(\mathcal{U}_{\infty})$, and more generally
$\pi_k(\mathcal{U}, \mathcal{U}_{\infty})\cong 0$ for all $k \leq
n$. Therefore, the same is true for any covering, in particular
for the universal abelian coverings: $\pi_k(\mathcal{U}^{ab},
\mathcal{U}_{\infty}^{ab})\cong 0$ for all $k \leq n$. Hence, by
Hurewicz Theorem, the vanishing also holds for the relative
homology groups, i.e., the maps of groups
$H_k(\mathcal{U}_{\infty}^{ab}) \to H_k(\mathcal{U}^{ab})$ are
isomorphism for $k<n$ and onto for $k=n$. Since these maps are
induced by an embedding (recall $n \geq 2$), the above are
morphisms of modules over the ring of Laurent polynomials in $s$
variables. The statement of the theorem follows now from
Lemma~\ref{div}.

\end{proof}

From now on to the end of this section, we will make the
assumption that the hyperplane at infinity $H$ is
\emph{transversal} (in the stratified sense) to the hypersurface
$V$. With this assumption, we show that the global cohomological
Alexander invariants of the hypersurface complement are entirely
determined by the degrees of the irreducible components on one
hand, and by the local topological information encoded by the
singularities of $V$ on the other hand. In particular, these
invariants depend on the local type of singularities of the
hypersurface.

First, we need some notations. Recall from \S~\ref{def} that
$A^q(\U) \cong H^q (\U, \LL^{\vee})$. For $x \in V$, we let
$\mathcal{U}_x=\mathcal{U} \cap B_x$, for $B_x$ a small open ball
at $x$ in $\mathbb{CP}^{n+1}$. Denote by $\mathcal{L}_x$ the
restriction of the local coefficient system $\mathcal{L}$ to
$\mathcal{U}_x$. Then the groups
$H^{\ast}(\mathcal{U}_x,\mathcal{L}^{\vee}_x)$ inherit a
$R_s$-module structure.

\begin{thm} \label{thm2}
Let $\lambda =(\lambda   _1,\cdots, \lambda  _s) \in \T^s$ and
$\epsilon \in \mathbb{Z}_{\geq 0}$. Fix an irreducible component
$V_1$ of $V$, and assume that $ \lambda    \notin
  \supp  (H^{q}(\mathcal{U}_x,\mathcal{L}^{\vee}_x))$ for all
$q<n+1-\epsilon$ and all points $x \in V_1$. Then $\lambda \notin
 \supp    (A^{q}(\mathcal{U}))$ for all $q<n+1-\epsilon$.
\end{thm}

\begin{proof}
Let $\mathcal{U}_1=\mathbb{CP}^{n+1} \setminus V_1$, and let
$i:\mathcal{U} \hookrightarrow \mathcal{U}_1$ and $j:
\mathcal{U}_1 \hookrightarrow \mathbb{CP}^{n+1}$ be the two
inclusions. Then $\mathcal{L}^{\vee}[n+1] \in
\text{Perv}(\mathcal{U})$, since $\mathcal{U}$ is smooth. Moreover
$\mathcal{F}:=Ri_{\ast}(\mathcal{L}^{\vee}[n+1]) \in
\text{Perv}(\mathcal{U}_1)$, since $i$ is a quasi-finite affine
morphism (see \cite{Sh}, Theorem 6.0.4). But $\mathcal{U}_1$ is
affine $(n+1)$-dimensional, and $\mathcal{F} \in
\text{Perv}(\mathcal{U}_1)$, therefore by Artin's vanishing
theorem for perverse sheaves (see \cite{Sh}, Corollary 6.0.4), the
following hold:

$$\mathbb{H}^k(\mathcal{U}_1, \mathcal{F})=0, \ \text{for all} \ k>0,$$
$$\mathbb{H}_c^k(\mathcal{U}_1, \mathcal{F})=0, \ \text{for all} \ k<0.$$
Let $a:\mathbb{CP}^{n+1} \to point$ be the constant map. Then:

$$\mathbb{H}^k(\mathcal{U}_1, \mathcal{F})\cong H^{k+n+1}(\mathcal{U},
\mathcal{L}^{\vee})\cong H^k(Ra_{\ast}Rj_{\ast}\mathcal{F})$$ and
$$\mathbb{H}_c^k(\mathcal{U}_1, \mathcal{F})\cong H^k(Ra_{!}Rj_{!}\mathcal{F})$$
Note that since $a$ is a proper map, we have $Ra_!=Ra_{\ast}$.

Now consider the canonical morphism $Rj_{!}\mathcal{F} \to
Rj_{\ast}\mathcal{F}$ and extend it to the distinguished triangle:
$$Rj_{!}\mathcal{F} \to Rj_{\ast}\mathcal{F} \to \mathcal{G}
\overset{[1]}{\to}$$ in $D_c^b(\mathbb{CP}^{n+1})$. Since
$j^{\ast}j_{!} \cong id \cong j^{\ast}j_{\ast}$, the complex
$\mathcal{G}$ is supported on $V_1$. Apply $Ra_!=Ra_{\ast}$ to the
above distinguished triangle and obtain:

$$Ra_{!}Rj_{!}\mathcal{F}
\to Ra_{\ast}Rj_{\ast}\mathcal{F} \to Ra_{\ast}\mathcal{G}
\overset{[1]}{\to}$$ Upon applying the cohomology functor to this
triangle, and using the above vanishing, we obtain that:

$$H^{k+n+1}(\mathcal{U},
\mathcal{L}^{\vee})\cong \mathbb{H}^k(\mathbb{CP}^{n+1},
\mathcal{G})\cong \mathbb{H}^k(V_1, \mathcal{G}) \ \ \text{for} \
\ k < -1,$$ and $H^{n}(\mathcal{U}, \mathcal{L}^{\vee})$ is a
sub-module of $\mathbb{H}^{-1}(V_1, \mathcal{G})$.

Therefore, by Lemma~\ref{div}, in order to prove the theorem it
suffices to show that, under our assumptions, $\lambda \notin
\text{Supp}(\mathbb{H}^{k}(V_1,\mathcal{G}))$ for all
$k<-\epsilon$. This  follows from the local calculation and the
hypercohomology spectral sequence. Indeed, for $x \in V_1$, we
have:
\begin{eqnarray*}
\mathcal{H}^q(\mathcal{G})_x &\cong &
\mathcal{H}^q(Rj_{\ast}\mathcal{F})_x \cong
\mathcal{H}^{q+n+1}(Rj_{\ast}Ri_{\ast}\mathcal{L}^{\vee})_x \cong
\mathbb{H}^{q+n+1}(B_x, R(j~\circ~i)_{\ast}\mathcal{L}^{\vee}) \\
&\cong & H^{q+n+1}(\mathcal{U}_x,\mathcal{L}^{\vee}_x)
\end{eqnarray*} where $\mathcal{U}_x=\mathcal{U} \cap B_x$, for
$B_x$ a small open ball at $x$ in $\mathbb{CP}^{n+1}$, and
$\mathcal{L}_x$ is the restriction of the local coefficient system
$\mathcal{L}$ to $\mathcal{U}_x$. Therefore, for a fixed $x \in
V_1$ the assumption that $ \lambda \notin
\text{Supp}(H^{q}(\mathcal{U}_x,\mathcal{L}^{\vee}_x))$ for all
$q<n+1-\epsilon$ is equivalent to the assumption $ \lambda  \notin
\text{Supp}(\mathcal{H}^{q}(\mathcal{G})_x)$ for all
$q<-\epsilon$. Next note that $\mathbb{H}^{k}(V_1,\mathcal{G})$ is
the abutment of a spectral sequence with the $E_2$-term defined by
$E_2^{p,q}=H^p(V_1,\mathcal{H}^q(\mathcal{G}))$. Moreover, if $
\lambda   \notin \text{Supp}(\mathcal{H}^{q}(\mathcal{G})_x)$ for
all $q<-\epsilon$ and for all $x \in V_1$, then $\lambda \notin
\text{Supp}(H^p(V_1,\mathcal{H}^{q}(\mathcal{G}))$ for
$p+q=k<-\epsilon$ (since $E_2^{p,q}$ is non-trivial only if $p
\geq 0$). Thus, from the spectral sequence, it follows that $
\lambda   \notin \text{Supp}(\mathbb{H}^{k}(V_1,\mathcal{G}))$ for
all $k<-\epsilon$. This finishes the proof of the theorem.

\end{proof}

\begin{remark}\label{rem2} \rm
In order to show that the universal cohomological modules depend
only on the local information around the singularities of the
hypersurface, it suffices to observe that the modules
$H^{*}(\mathcal{U}_x,\mathcal{L}^{\vee}_x)$, $x \in V_1$, are
entirely determined by the local universal homological Alexander
modules at $x$.\newline In order to see this, we first introduce
some notation: let $\U_0$ denote the hypersurface complement
$\mathbb{CP}^{n+1} \setminus V$, and for $x \in V_1$ we set
$\U'_x=\U_0 \cap B_x$, for $B_x$ a small open ball at $x$ in
$\mathbb{CP}^{n+1}$. Note that $H_1 (\U'_x)=\Z^k$, where $k$ is
the number of irreducible components of the hypersurface
singularity germ $(V,x)$ (cf. \cite{Di}, p.103). Let $\U_x ^{ab}$
and $(\U'_x) ^{ab}$ be the universal abelian covers of $\U_x$ and
$\U'_x$, respectively, and denote by $A_{\ast}(\U_x)$ and
respectively $A_{\ast}(\U'_x)$ the associated universal
homological Alexander modules. The modules $A_{\ast}(\U'_x)$ will
be called the local universal homological Alexander modules at
$x$, as they depend only on the singularity germ $(V,x)$.

If $i_x : \U_x \hookrightarrow \U$ denotes the inclusion map, then
the local system $\LL_x$ on $\U_x$ is induced via the composition
of maps
$$\phi : \pi_1(\U_x) \overset{(i_x)_{\#}}{\to} \pi_1(\U)
\overset{\text{lk}}{\to} H_1(\U) \to \text{Aut}(R_s)$$ On the
other hand, by the naturality of the Hurewicz morphism, $\phi$
factors through $\text{lk}_x :\pi_1(\U_x) \to H_1(\U_x)$, $R_s$
becoming in this way a (left) $\C[H_1 (\U_x)]$-module. Then, by
\cite{Di2} p.50, it follows that
$H^{*}(\mathcal{U}_x,\mathcal{L}^{\vee}_x)$ is the homology of the
equivariant Hom:

$$C^{*}(\mathcal{U}_x,\mathcal{L}^{\vee}_x)=\text{Hom}_{\C[H_1 (\U_x)]}(C_{\ast}^0(\U_x ^{ab}), R_s),$$
where $R_s$ is regarded now as a right $\C[H_1 (\U_x)]$-module
using the involution on the group ring as in Remark \ref{rem0},
and as a left $R_s$-module. By \cite{L}, p.6, there is a spectral
sequence converging to $H^{*}(\mathcal{U}_x,\mathcal{L}^{\vee}_x)$
with

\begin{equation} \label{spsq10}
E_2^{p,q}=\text{Ext}^q _{\C[H_1 (\U_x)]} (A_p(\U_x),R_s).
\end{equation}

In order to fully justify our claim, it remains to relate the
local universal Alexander invariants $A_{\ast}(\U'_x)$ to the
modules $A_{\ast}(\U_x)$, at points $x \in V_1$.\newline For
points $x \in V_1 \setminus (V_1 \cap H)$ we have $\U'_x=\U_x$,
thus our claim follows for such points by the above spectral
sequence.\newline If $x \in V_1 \cap H$ then due to the
transversality assumption, it's easy to see that $\U_x$ is
homotopy equivalent to $\U'_x \times S^1$. It follows that
$\U_x^{ab} \simeq (\U'_x)^{ab} \times \mathbb{R}$, thus by the
homological K\"{u}nneth formula we obtain that the group
$A_p(\U_x)$ is isomorphic to $H_p((\U'_x)^{ab},\C) \otimes
H_0(\mathbb{R}, \C) \cong A_p(\U'_x)$. When regarded as a $\C[H_1
(\U_x)]$-module, the isomorphism can be written as (see \cite{CS},
Prop. 1.8):

$$A_p(\U_x) \cong (A_p(\U'_x) \otimes_{\C[H_1 (\U'_x)]} \C[H_1 (\U_x)]) \otimes_{\C[\Z]} \C.$$
Together with the spectral sequence \ref{spsq10} this finishes the
proof of the claim.

\end{remark}

\begin{remark}\rm

If $S$ is an $s$-dimensional stratum in a Whitney stratification
of $V$ such that $x \in S$, then $A_{p}(\U'_x)=0$ if $p > n-s$.
Indeed, $\U'_x$ has the homotopy type of the link complement
$S^{2n-2s+1}_x \setminus L_x$, where $S^{2n-2s+1}_x$ is a small
sphere at $x$ in a submanifold of $\mathbb{CP}^{n+1}$ which meets
$S$ transversally at $x$ (and no other point), and $(S^{2n-2s+1}_x
, L_x)$ is the link pair of  the stratum $S$ in the pair
$(\mathbb{CP}^{n+1},V)$. Since $S^{2n-2s+1}_x \setminus L_x$
admits a cyclic cover which has the homotopy type of a CW complex
of dimension $n-s$ (i.e. the fiber of the Milnor fibration
associated to the algebraic link $(S^{2n-2s+1}_x , L_x)$), it
follows that the universal abelian cover $(\U'_x)^{ab}$ has the
homotopy type of a $(n-s)$-dimensional CW complex, thus proving
the claim.

\end{remark}

\bigskip

The following consequence of Theorem~\ref{thm2},
Remark~\ref{rem2}, and of Example~\ref{ex1} is similar to some
results in \cite{DiLi}, \cite{Li8}, \cite{Li9}.

\begin{cor} \label{cor1}
\ \\
\noindent (i) (Case $\epsilon=0$) With the notation in the above theorem, assume in addition that
$V$ is a normal crossing divisor at any point of the component
$V_1$. Then $\supp (A^k(\U)) \subset \{ {\bf 1} \}$ for any
$k<n+1.$

\medskip

\noindent (ii) (Case $\epsilon=1$) With the notation in the above theorem, assume in addition that
$V$ is an INNC divisor at any point of the component
$V_1$. Then $\supp (A^k(\U)) \subset \{ {\bf 1} \}$ for any
$k<n.$

\end{cor}

Using a similar argument (see also \cite{Di4}) we obtain the
following result.

\begin{thm}\label{ucam}
Assume that the hypersurface $V$ is transversal (in the stratified
sense) to the hyperplane at infinity $H$. Then for $k \leq n$,
$  \supp   (A^k(\mathcal{U}))$ is contained in the zero set of
the polynomial $t_1^{d_1}\cdots t_s^{d_s}-1$, thus has positive
codimension in $\T^s$.
\end{thm}

The positive codimension property of supports in the
universal abelian case should be regarded as the analogue of the
torsion property in the infinite cyclic case (cf. \cite{Max},
\cite{Di4}). Example \ref{ex3} below shows that transversality except
finitely many points is not enough to get Theorem \ref{ucam}.

\begin{proof}
As in the proof of the previous theorem, after replacing
$\mathcal{U}_1$ by the affine space
$\mathbb{C}^{n+1}=\mathbb{CP}^{n+1} \setminus H$, it follows that
for $k \leq -1$, $H^{k+n+1}(\mathcal{U}, \mathcal{L}^{\vee})$ is a
sub-module of $\mathbb{H}^{k}(\mathbb{CP}^{n+1}, \mathcal{G})$,
where $\mathcal{G}$ is now a complex of sheaves supported on $H$.
Therefore, by Lemma~\ref{div}, it suffices to prove the theorem
for the supports of the modules $\mathbb{H}^{k}(H, \mathcal{G})$
with $k \leq -1$.

As in the previous theorem, for $x \in H$, the local calculation
on stalks yields $\mathcal{H}^q(\mathcal{G})_x \cong
H^{q+n+1}(\mathcal{U}_x,\mathcal{L}^{\vee}_x)$, where
$\mathcal{U}_x=\mathcal{U} \cap B_x$, for $B_x$ a small open ball
at $x$ in $\mathbb{CP}^{n+1}$. If $x \in H \setminus H\cap V$,
then $\mathcal{U}_x$ is homotopy equivalent to
$\mathbb{C}^{\ast}$, and the corresponding local system
$\mathcal{L}^{\vee}_x$ is defined by the action of
$\gamma_{\infty}$, i.e. by  multiplication by $\prod_{j=1}^s
{(t_j)^{d_j}}$. On the other hand, if $x \in V \cap H$, then due
to the transversality assumption, $\mathcal{U}_x$ is homotopy
equivalent to a product $(B'_x \setminus V \cap B'_x) \times
\mathbb{C}^{\ast}$, with $B'_x$ a small open ball centered at $x$
in $H$, and the local system $\mathcal{L}^{\vee}_x$ is an external
tensor product, the second factor being defined by the
multiplication by $\prod_{j=1}^s {(t_j)^{d_j}}$. Thus, by the
Kunneth spectral sequence, the stalk cohomology groups of
$\mathcal{G}$ along $H$, i.e. $\mathcal{H}^q(\mathcal{G})_{x\in
H}$, have supports contained in the zero set of the polynomial
$t_1^{d_1}\cdots t_s^{d_s}-1$. Then by the hypercohomology
spectral sequence, the same is true for the supports of the
hypercohomology groups $\mathbb{H}^{k}(H, \mathcal{G})$.

\end{proof}

\section{Explicit Computations via Logarithmic Connections}

We review a general method used to determine the characteristic varieties in the case of hyperplane arrangements, see
\cite{ESV} and \cite{STV},  and show that essentially the same method applies to more general situations as well.

Let $\pi : (Z,D) \to (\mathbb{CP}^{n+1},V \cup H)$ be an embedded resolution of singularities for the reduced divisor
$V \cup H$. In particular

\noindent (i) $D$ is a normal crossing divisor with smooth irreducible components;

\noindent (ii) $\pi : Z \setminus D \to \U$ is an isomorphism.

\noindent In this setting there is a Hodge-Deligne spectral sequence
\begin{equation}\label{Deligne1}
E_1^{p,q}=H^q(Z, \Omega^p _Z(logD)) \Rightarrow H^{p+q}(\U,\C)
\end{equation}
degenerating at $E_1$ and inducing the Hodge filtration $F$ of the Deligne mixed Hodge structure on
$H^{p+q}(\U,\C)$, see \cite{De2}.

\begin{cor}\label{Delcor}

If the Deligne mixed Hodge structure on some cohomology space $H^{m}(\U)$ is pure of type $(m,m)$, then

\noindent (i) $H^0(Z, \Omega^m _Z(logD))=H^{m}(\U)$ and

\noindent (ii) $H^q(Z, \Omega^p _Z(logD))=0$ for $p+q=m$ and $q>0$.

\end{cor}

We list below several cases when this property holds.

\begin{example}\label{expure1} \rm

\noindent (a)  When $V$ is a hyperplane arrangement, the cohomology space $H^{m}(\U)$ is pure of type $(m,m)$     for all $m \geq 0$, see  \cite{DiL}.

\noindent (b)  When $V$ is a smooth rational curve arrangement in the projective plane (i.e. any irreducible component of $V$ is either a line or a smooth conic), the cohomology space $H^{m}(\U)$ is pure of type $(m,m)$  for all
$m \geq 0$ (easy exercise for the reader).

\noindent (c) $H^m(\U)$ is always pure of type $(m,m)$ for all $m \leq 1$. This follows from the fact that
$g=(g_1,...,g_s):\U \to \T^s$ induces an isomorphism at the $H^m$-level all $m \leq 1$.
Here we look at $\U$ as a subset of $\C^{n+1}$ and we set $g_j(x_1,...,x_{n+1})=f_j(1,x_1,...,x_{n+1})$.

\end{example}

For $\lambda =( \lambda _1,..., \lambda _s) \in \T^s$, let $\LL_{\lambda}$ be the corresponding local system on $\U= Z \setminus D  $. Let $\alpha _j \in \C$ be such that $\exp(-2\pi i \alpha _j)= \lambda _j$ for $j=1,...,s.$
Then $\LL_{\lambda}$ is the local system of horizontal sections of the connection
$$ \nabla _{\alpha}: \OO _{\U} \to \Omega ^1_{\U}$$
given by $ \nabla _{\alpha}(u)=du+ u \cdot \omega _{\alpha}$ where
$$ \omega _{\alpha}=\sum _{j=1,s} \alpha _j \frac{dg_j}{g_j}.$$
 Alternatively,
if we look at $\U$ as a subset of $\mathbb{CP}^{n+1}$, then we can use the formula
$$ \omega _{ \alpha  }=\sum _{j=0,s} \alpha _j \frac{df_j}{f_j}$$
where we set $ \alpha _0=- \sum _{j=1,s} d_j \cdot \alpha _j$. Recall that $f_0=x_0$.

\medskip

Using the fact that $\U$ is affine and our connection is regular, it follows that
\begin{equation} \label{deRham1}
H^m(\U,\LL_{\lambda})=H^m(H^0(\U, \Omega^*  _{\U}), \nabla _{\alpha} )
\end{equation}
just as in  \cite{Di2}, (Thm. 3.4.18) or, for complete proofs, \cite{De1}. However, this result is not so useful
to perform explicit computations since the groups $H^0(\U, \Omega^*  _{\U})$ are too large.

There is a second approach to computing $H^m(\U,\LL_{\lambda})$,
this time using {\it logarithmic connections}. It has the
advantage of reducing the size of the spaces  $H^0(\U, \Omega^*
_{\U})$, but one has to be more careful about the {\it residues}
$\alpha _j$. More precisely, the pull-back of the connection $
\nabla _{ \alpha }$ under the embedded resolution $\pi$ is a
logarithmic connection $\tilde  \nabla _{\alpha  }$ on $Z$ with
poles along $D$. Let $\rho_i$ be the residue of the connection
$\tilde  \nabla _{ \alpha }$ along the irreducible component $D_i$
of $D$. When $D_i$ is the proper transform of some component $V_j$
of $V$ one has $\rho_i=\alpha _j$.
 
\begin{defn} \label{defadmis}

A choice of residues  $\alpha=( \alpha_0,
\alpha_1,...,\alpha_s)$ for $\LL_{\lambda}$ as above  is an {\it admissible choice of residues}
for $\LL_{\lambda}$ if $\rho_i \notin \N_{>0}$ for all irreducible
components $D_i$ of $D$. A rank one local system $\LL_{\lambda}$ is  {\it admissible} if
there is some admissible choice of residues for it.

\end{defn}

\begin{remark} \label{admis1rk} \rm

It is easy to see using Hironaka's embedded resolution of singularities by blowing-up smooth subvarieties,
that for any $i$ there is a relation
$$\rho_i=\sum_{j=1,s}n_{ij}\alpha_j$$
with $n_{ij} \in \Z$ (see \cite{ESV} for similar formulas and note that negative coefficients occur due to the presence of the hyperplane at infinity). The condition $\rho_i \notin \N_{>0}$ is clearly satisfied if all $\alpha_j$
are sufficiently small. In other words, there is a neighborhood $U({\bf 1})$ of the trivial local system
${\bf 1} \in \T^s$ formed entirely by admissible local systems.

If we move away from the trivial local system , it is not clear whether all the local systems are admissible.
 The answer to this
question is negative for some hyperplane arrangements, see
\cite{CS1}, Example 4.4, \cite{CO}, Example 3.4 and  \cite{LY}. On the other hand, for not very complicated arrangements, see Examples \ref{linesandconic1} and 
 \ref{linesandconic2} below, the answer is positive.

\end{remark}

 For an  admissible choice of residues one
has an $E_1$-spectral sequence
\begin{equation} \label{Deligne2}
E_1^{p,q}=H^q(Z, \Omega^p _Z(logD)) \Rightarrow H^{p+q}(\U, \LL_{\lambda}   )
\end{equation}
whose differential $d_1$ is induced by $\tilde  \nabla _{ \alpha  }$, see  \cite{Di2}, (Thm. 3.4.11 (i)). The above discussion proves the following.

 \begin{prop}\label{proppure}
Assume that  $\alpha=( \alpha_0, \alpha_1,...,\alpha_s)$ is an admissible choice of residues for $\LL_{\lambda}$ and that the cohomology groups $H^m(\U)$ are pure of type $(m,m)$ for all
$m \leq k$. Then
$$ H^m(\U, \LL_{\lambda}   )=H^m(H^*(\U),  \omega _{\alpha }    \wedge )$$
for all $m \leq k$ and $H^{k+1}(H^*(\U),  \omega _{ \alpha }    \wedge )$ is a subspace in  $H^{k+1}(\U, \LL_{\lambda}   )$.

\end{prop}

When $\U$ is a hyperplane arrangement complement, this is exactly the argument used in \cite{ESV} and \cite{STV}.
Proposition \ref{proppure}, Remark \ref{admis1rk}  and Example \ref{expure1} yield the following.

\begin{cor}\label{1charvar}
If $\U$ is any affine hypersurface arrangement complement, then there is a neighborhood 
 $U({\bf 1})$ of the trivial local system
${\bf 1} \in \T^s$ such that 
$$ H^1(\U, \LL_{\lambda}   )=H^1(H^*(\U),  \omega _{\alpha }    \wedge )$$
for  any local system $\LL_{\lambda} \in U({\bf 1})$, $\alpha$ being an arbitrary choice of admissible residues for $\LL_{\lambda}$.

\end{cor}

\begin{cor}\label{hyperplanes}
If $\U=M(\A)$ is a hyperplane arrangement complement, then there is a neighborhood 
 $U({\bf 1})$ of the trivial local system
${\bf 1} \in \T^s$ such that 
$$ H^m(\U, \LL_{\lambda}   )=H^m(H^*(\U),  \omega _{\alpha }    \wedge )$$
for any $m \in \N$, and any local system $\LL_{\lambda} \in U({\bf 1})$, $\alpha$ being an arbitrary choice of admissible residues for $\LL_{\lambda}$.

\end{cor}

 \begin{example} \label{linesandconic1} \rm
In the projective plane $\mathbb{CP}^{2}$ consider the hypersurface $V$ having as irreducible components $V_1:x=0$, $V_2: y=0$, $V_3:x^2-yz=0$. Let
$H=V_0$ be the line at infinity given by $z=0$ and note that $H$ is {\it not transverse} in a stratified sense to $V$.
Consider the connection  $ \nabla _{\lambda}$ whose residues are
 $\alpha=( \alpha_0, \alpha_1,...,\alpha_3)$ with
$$ \alpha_0=-\alpha_1 -\alpha_2-2 \alpha_3.$$
Let $A=V_1 \cap V_2 \cap V_3=(0:0:1)$ and $B=V_1 \cap V_0 \cap V_3=(0:1:0)$.
To construct the embedded resolution of $V \cup H$ we first blow-up the points $A$ and $B$, creating thus two exceptional divisors, $D_A$ and respectively
$D_B$. The corresponding residues along  $D_A$ and $D_B$ are easily computable and we get $ \alpha_A=\alpha_1 +\alpha_2+\alpha_3$
 and respectively   $ \alpha_B=\alpha_1 +\alpha_0+\alpha_3= -\alpha_2- \alpha_3. $
Let $P=D_A  \cap V_2' \cap V_3 '$ and $Q=D_B\cap V_0' \cap V_3'$, where $'$ denotes the proper transform of a divisor.
To get the embedded resolution of $V \cup H$ we have just to  blow-up the points $P$ and $Q$, creating thus two new exceptional divisors, $D_P$ and respectively
$D_Q$.  The corresponding residues are $ \alpha_P= -\alpha_Q=\alpha_1 +2\alpha_2+2\alpha_3$.
Therefore the choice of residues  $\alpha=( \alpha_0, \alpha_1,...,\alpha_3)$
is admissible if and only if none of the residues
$$ \alpha_1,  \alpha_2,  \alpha_3, -\alpha_1 -\alpha_2-2 \alpha_3, \alpha_1 +\alpha_2+\alpha_3,-\alpha_2- \alpha_3,\alpha_1 +2\alpha_2+2\alpha_3, -(   \alpha_1 +2\alpha_2+2\alpha_3)$$
is a strictly positive integer.

\begin{lemma} \label{admis1}
In the situation of Example \ref{linesandconic1}, any rank one local system is admissible.

\end{lemma}

\proof

It is clearly enough to consider the case of real residues $\alpha_j$. Otherwise, we just look at the corresponding real parts.

 We divide the possibilities in the following two cases.

\noindent Case 1. ($\alpha_1 +2\alpha_2+2\alpha_3 \notin \Z$).

Suppose first that in addition $\alpha_1 +\alpha_2+\alpha_3 \notin \Z$. Then the choice with  $\alpha_j \in [0,1)$ for $j=1,2,3$ is admissible.

Now suppose that  $\alpha_1 +\alpha_2+\alpha_3 \in \Z$. It follows that $\alpha_2+\alpha_3 \notin \Z$.
 Then the choice with  $\alpha_j \in [0,1)$ for $j=2,3$ and $ \alpha_1<0$ such that  $\alpha_1 +\alpha_2+\alpha_3 =0$   is admissible.

\noindent Case 2. ($\alpha_1 +2\alpha_2+2\alpha_3 \in \Z$).

Then we have to choose $\alpha_1= -2\alpha_2-2\alpha_3$. The residues in this case are just
$$-2( \alpha_2+\alpha_3), ~~ -( \alpha_2+\alpha_3), ~~ \alpha_2, ~~\alpha_3.$$
Hence it is enough to take 
 $\alpha_j \in [0,1)$ for $j=2,3$.

\endproof

Now we continue Example  \ref{linesandconic1} by applying Example \ref{expure1} and  Proposition \ref{proppure} to get $ H^m(\U, \LL_{\lambda}   )=H^m(H^*(\U),  \omega _{\alpha  }    \wedge )$
for all $m$. In order to perform this computation, we need a precise description of the cohomology algebra $ H^*(\U)$ (with $\C$ coefficients) and this can be obtained in this example from the local considerations in \cite{Di}, pp. 47-49.
The result can be described as follows.

\noindent (i) $H^0(\U)=\C$ and the generator is $1$;

\noindent (ii) $H^1(\U)=\C^3$ and a basis is given by $\eta _1=\frac{dx}{x}$,
$\eta _2=\frac{dy}{y}$ and $\eta _3=\frac{d(x^2-y)}{x^2-y}$;

\noindent (iii) $H^2(\U)=\C^2$ and a basis is given by
$\eta _{12}=\eta _1  \wedge \eta _2  $ and $\eta _{23}=\eta _2  \wedge \eta _3  $. The multiplication is given by the relation
$$2\eta _1  \wedge \eta _3=2\eta _{12}+ \eta _{23}.$$

\noindent (iv) Since $\U$ is affine, $H^m(\U)=0$ for $m>2$.

The computation of $H^m(H^*(\U),  \omega _{\lambda}    \wedge )$ falls into 3 cases.

\medskip

\noindent Case 1. $\alpha_1=  \alpha_2=  \alpha_3=0$ and $ \LL_{\lambda}=\C$ is the constant local system. Then of course  $ H^m(\U, \LL_{\lambda}   )=H^m(\U)$ for all $m$.

\medskip

\noindent Case 2. $\alpha_1 +2\alpha_2+2\alpha_3=0$. Then a direct computation shows that  $ H^0(\U, \LL_{\lambda}   )=0$ and
$\dim H^1(\U, \LL_{\lambda}   )=\dim H^2(\U, \LL_{\lambda}   )=1.$

\medskip

\noindent Case 3. $\alpha_1 +2\alpha_2+2\alpha_3\ne 0$. Again a direct computation shows that
$ H^0(\U, \LL_{\lambda}   )=   H^1(\U, \LL_{\lambda}   )= H^2(\U, \LL_{\lambda}   )=0.$

\medskip

The above computations yield  the
following equalities. 

$$V_t^{0,1}(\U)= V_t^{0,2}(\U)   =\{ {\lambda} \in \T^3; \  \ \lambda_1 \lambda_2^2 \lambda_3^2=1 \}$$
$$V_t^{1,1}(\U)=V_t^{2,1}(\U)=V_t^{1,2}(\U)=\{ {\bf 1} \}$$
$V_t^{m,1}(\U)= \emptyset$ for $m>2$ and $V_t^{m,2}(\U)= \emptyset$ for $m>1$.

These results are consistent with the general results by Arapura \cite{A}. See also Suciu
 \cite{Su} for a related discussion.

 Note that the  above 2-dimensional subtorus $\T=\{ {\lambda} \in \T^3; \  \ \lambda_1 \lambda_2^2 \lambda_3^2=1 \}$     is different from the 2-dimensional subtorus predicted by Theorem \ref{ucam} in the case of a divisor $V$ transversal to the line at infinity.

A special class of local systems is formed by the equimonodromical local systems $ \LL_{\lambda} $
such that $\lambda_0=\lambda_1= ... =\lambda_3$. Then, for $\lambda_0^5=1$, the dimension of the cohomology space $H^m(\U, \LL_{\lambda}   )$ is exactly the multiplicity of the root $t=\lambda_0$ in the
characteristic polynomial
$$\Delta ^m(t)=det(t\cdot Id - h^m)$$
where $F:xyz(x^2-yz)=1$ is the associated Milnor fiber in $\C^3$ and $h:F \to F$ is the monodromy operator, see for instance \cite{Di2}, (6.4.6). To compute the cohomology of such a  equimonodromical local system $ \LL_{\lambda} $, one should start by an admissible choice for the residues $\alpha=(\alpha_0, \alpha_1,\alpha_2, \alpha_3)   $. For instance, the obvious choice $\alpha = ( \frac{-4}{5}, \frac{1}{5}, \frac{1}{5},   \frac{1}{5}) $ is not admissible. A good choice here is  $\alpha = ( \frac{1}{5}, \frac{-4}{5},  \frac{1}{5},    \frac{1}{5}) $. Using this choice, we get the following characteristic polynomials
in this situation.

$$\Delta ^0(t)=t-1, ~ ~ ~ \Delta ^1(t)=(t-1)^2(t^5-1), ~ ~ ~ \Delta ^2(t)=(t-1)(t^5-1).$$

\end{example}

The following example is similar to the previous one, but it exhibits a curve $V$ which is transversal to the line at infinity $H$ and it needs a different approach for the computation of the cohomology algebra $H^*(\U)$. Moreover, in this case the cohomology algebra $H^*(\U)$ is not spanned by the degree one part $H^1(\U)$.

\begin{example} \label{linesandconic2} \rm

In the projective plane $\mathbb{CP}^{2}$ consider the hypersurface $V$ having as irreducible components $V_1:x=0$, $V_2: y=0$, $V_3:x^2-y^2+yz=0$. Let
$H=V_0$ be the line at infinity given by $z=0$ and note that $H$ is {\it transverse} in a stratified sense to $V$.
Consider the connection  $ \nabla _{\lambda}$ whose residues are
 $\alpha=( \alpha_0, \alpha_1,...,\alpha_3)$ with
$$ \alpha_0=-\alpha_1 -\alpha_2-2 \alpha_3.$$
Let $A=V_1 \cap V_2 \cap V_3=(0:0:1)$.
To construct the embedded resolution of $V \cup H$ we first blow-up the point $A$, creating an exceptional divisor $D_A$. The corresponding residue along  $D_A$ is $ \alpha_A=\alpha_1 +\alpha_2+\alpha_3$.
Let $P=D_A  \cap V_2' \cap V_3 '$, where $'$ denotes the proper transform of a divisor.
To get the embedded resolution of $V \cup H$ we have just to  blow-up the point $P$, creating a new exceptional divisor $D_P$.  The corresponding residue is $ \alpha_P=\alpha_1 +2\alpha_2+2\alpha_3$.
Therefore the choice of residues  $\alpha=( \alpha_0, \alpha_1,...,\alpha_3)$
is admissible in this case if and only if none of the residues
$$ \alpha_1,  \alpha_2,  \alpha_3, -\alpha_1 -\alpha_2-2 \alpha_3, \alpha_1 +\alpha_2+\alpha_3,\alpha_1 +2\alpha_2+2\alpha_3$$
is a strictly positive integer. It can be shown, exactly as in Lemma \ref{admis1} above, that in this situation any rank one local system is admissible.

It follows that we can apply Example \ref{expure1} and  Proposition \ref{proppure} to get $ H^m(\U, \LL_{\lambda}   )=H^m(H^*(\U),  \omega _{\lambda}    \wedge )$
for all $m$. To get a precise description of the cohomology algebra $ H^*(\U)$ we can proceed as follows.

\noindent (i) $H^0(\U)=\C$ and the generator is $1$;

\noindent (ii) $H^1(\U)=\C^3$ and a basis is given by $\eta _1=\frac{dx}{x}$,
$\eta _2=\frac{dy}{y}$ and $\eta _3=\frac{d(x^2-y^2+y)}{x^2-y^2+y}$;

\noindent (iii) To compute $H^2(\U)$ is the first difficulty. This can be done by setting
$\U^0 = \mathbb{CP}^{2} \setminus (V_0 \cup V_1 \cup V_2)$, $V_3^0=V_3 \setminus (V_0 \cup V_1 \cup V_2)$ and considering the Gysin sequence
$$ H^1(\U) \to  H^0(V_3^0   ) \to H^2(\U^0) \to  H^2(\U) \to H^1(V_3^0   ) \to 0.$$
The first morphism, given by the Poincar\'e-Leray residue $R$, is clearly surjective, i.e.
$R(\eta _3)=1$. Then $\dim  H^2(\U^0)=1$ and a generator is $\eta_{12}=\eta _1 \wedge \eta _2$.
The affine curve $V_3^0$ is isomorphic to $\C \setminus \{-1,0,1\}$ under the parametrization
$$ x= \frac{t}{t^2-1}, ~  ~   y= \frac{t^2}{t^2-1}.$$
Using this parametrization, we can identify $H^1(V_3^0   )$ to $\C^3$ by sending a rational differential form to its residues at the points $ \{-1,0,1\}$. Some explicit computations involving the last nonzero morphism in the exact sequence above (which is again given by the Poincar\'e-Leray residue $R$)  show that
$R(\eta_{13})$ and $R(\eta_{23})$ are linearly independent in $H^1(V_3^0   )=\C^3$, where $\eta _{13}=\eta _1  \wedge \eta _3  $ and $\eta _{23}=\eta _2  \wedge \eta _3  $.
It follows that $\eta _{12}$, $\eta_{13}$ and $\eta _{23}$ are linearly independent in $H^2(\U)$,
which is 4-dimensional.

It follows that the following cases are possible in this example.

\medskip

\noindent Case 1. $\alpha_1=  \alpha_2=  \alpha_3=0$ and $ \LL_{\lambda}=\C$ is the constant local system. Then of course  $ H^m(\U, \LL_{\lambda}   )=H^m(\U)$ for all $m$.

\medskip

\noindent Case 2. $(\alpha_1,\alpha_2,\alpha_3)\ne (0,0,0)$. Then a direct computation shows that  $ H^0(\U, \LL_{\lambda}   )= H^1(\U, \LL_{\lambda}   ) =0$ and
$\dim H^2(\U, \LL_{\lambda}   )=2.$

\medskip

The above computations yield  the following equalities.

$$ V_t^{0,2}(\U)=  V_t^{1,2}(\U)   =\T^3$$
(hence here the support has 0-codimension)
$$V_t^{1,1}(\U)=V_t^{2,1}(\U)=V_t^{2,2}(\U) = V_t^{3,2}(\U)  =\{ {\bf 1} \}$$
$V_t^{m,1}(\U)= \emptyset$ for $m>2$ and $V_t^{m,2}(\U)= \emptyset$ for $m>3$.
Note that the inclusion in Theorem \ref{ucam} is strict in this case.

Consider as in the above example  the associated Milnor fiber   $F:xyz (x^2-y^2+y)   =1$  and $h:F \to F$  the monodromy operator.  A good choice of residue is again given by   $\alpha = ( \frac{1}{5}, \frac{-4}{5},  \frac{1}{5},    \frac{1}{5}) $. Using this choice, we get the following characteristic polynomials
in this situation.

$$\Delta ^0(t)=t-1, ~ ~ ~ \Delta ^1(t)=(t-1)^3, ~ ~ ~ \Delta ^2(t)=(t-1)^2(t^5-1)^2.$$

\end{example}

\begin{remark} \label{localglobal} \rm
In order to apply Theorem \ref{thm2}, we have to check the
vanishing of some local cohomology groups. When the hypersurface
germs occuring in these local complements are quasi-homogeneous,
then we can globalise the local situation and compute the
corresponding  local cohomology groups using the ideas explained
in this section. For instance, Example  \ref{linesandconic1}
covers the case of a plane curve singularity consisting of 3
smooth branches $(C_1,0)$, $(C_2,0)$ and $(C_3,0)$ such that the
intersection multiplicities are given by $(C_1,C_2)=1$,
$(C_1,C_3)=1$ and $(C_2,C_3)=2$. This follows from the topological
classification of the plane curve germs, see
 \cite{Di}, p. 45.

\end{remark}

\section{A More General Setting}
In this section we define multi-variable Alexander invariants in a
more general setting (see below) and attempt to relate them to the
invariants previously defined.

Assume that the  hypersurface  $V$ in $\mathbb{CP}^{n+1}$ has $s$
irreducible components $V_i$  with degrees $\text{deg}(V_i)=d_i$
for $i=1,\cdots,s$. Denote by $\U _0$ the complement
$\mathbb{CP}^{n+1} \setminus V$, and let $d=g.c.d.(d_1,\cdots,d_s)$. Then
$$H_1(\U _0) = \Z^{s-1} \oplus (\Z / d\Z)$$ is generated by the
meridians $\gamma_i$ about the non-singular part of each component
$V_i$, for $i=1,\cdots,s$ (cf. \cite{Di}, (4.1.3)). These
meridians satisfy a single relation, namely $$\sum_{i=1}^s {d_i
\gamma_i} = 0.$$ Now fix a hyperplane $H$ and set, as before, $\U
= \mathbb{CP}^{n+1} \setminus (V \cup H)$. Recall that $H_1 (\U) = \Z^s$,
freely generated by the meridians $\gamma_i$, $i=1,\cdots,s$. Let
$i: \U \hookrightarrow \U_0$ be the inclusion map, and denote by
$\U_0 ^{ab}$ and $\U^{ab}$ the universal abelian covers of $\U_0$
and $\U$ respectively, and by $p_0$ and $p$ the corresponding
covering projections.

The invariants we are interested in are $H_\ast (\U_0 ^{ab} ;
\C)$, regarded as modules over the quotient ring $\C[H_1(\U_0)] =
\C[t_1 ^{\pm1},\cdots, t_s ^{\pm1}]/(t_1^{d_1}\cdots
t_s^{d_s}-1)$.

\bigskip

It is a natural question to  find the relation between the universal abelian
invariants associated to the complement of $V$, and those
associated to the complement of $V \cup H$.

\bigskip

For a topological space $X$, let $\LL (X)$ denote the set of rank one complex local systems on $X$. When $X=\U$, then  $\LL (\U)$ is naturally identified to the $s$-dimensional complex torus $\T^s$. For  $X=\U _0$, the set  $\LL (\U _0)$
corresponds to the subset in $\T^s$ given by
$$ \{ \lambda=( \lambda _1,..., \lambda_s) \in \T^s~~ | \lambda _1^{d_1}\cdots
 \lambda _s^{d_s}=1\}.$$
With the notation above, let $d_j=d \cdot d_j'$ and consider the $(s-1)$-dimensional complex subtorus
$$\T = \{ \lambda=( \lambda _1,..., \lambda_s) \in \T^s~~ | \lambda _1^{d_1'}\cdots
 \lambda _s^{d_s'}=1\}.$$
For each $d$-root of unity $\beta$, let $\lambda (\beta  )$ be one point in the
hypersurface in $\T^s$ given by the equation
$$ \lambda _1^{d_1'}\cdots
 \lambda _s^{d_s'}=\beta  .$$
Then  $\LL (\U _0)$ is precisely the disjoint union of translated tori given by
$$\LL (\U _0) = \bigcup _{\beta  } \lambda (\beta  ) \T.$$
A different way of looking at a local system $\LL$ in  $\LL (\U _0)$ is by considering it as a local system in  $\LL (\U )$ (given by the obvious restriction $\LL | \U$) such that the action of the elementary loop about the hyperplane $H$ is trivial. This view-point yields the following exact sequence
$$ \cdots \to H^k(\U _0,\LL)  \to H^k(\U ,\LL) \to H^{k-1}(\U _0 \cap H ,\LL)
\to H^{k+1}(\U _0,\LL)\to  \cdots$$
for details on this see \cite{Di2}, pp. 221-222. The following obvious consequence should be compared to
\cite{Ra}, \cite{LV}, Proposition 1.3. The higher dimensional case, but with a generic hyperplane at infinity $H$, was considered in \cite{Li2}, Lemmas 1.5, 1.11 and 1.13.

\begin{cor} \label{curvearr}
Assume that $V$ is a plane curve arrangement, i.e. $n=1$. Then, for any rank one local system
$\LL= \LL_{ \lambda}  $ on $\U_0$ and any choice of the line at infinity $H$, one has
$$ \dim H^1(\U,\LL)= \dim H^1(\U_0,\LL)+ \epsilon.$$
Here $\epsilon \in \{0,1\}$ and $\epsilon=0$ if and only if there is a point $p\in V\cap H$ such that
$$\prod_{j=1,s} \lambda _j^{k_j} \ne 1$$
where $k_j=mult_p(V_j,H)$ is the intersection multiplicity of the component $V_j$ and the line $H$ at the point $p$.

\end{cor}

One case which is already well-explored is the following.

\begin{example} \label{onecomp} \rm
Assume that $n>1$, $s=1$ and that $V=V_1$ is a hypersurface of degree $d$ having only isolated singularities.
Then $\pi _1(\U_0)=\Z/d\Z$ and hence a local system $\LL= \LL_{ \beta} $ corresponds to a choice of a $d$-root of unity $\beta$. For $\beta =1$ we get $H^0(\U_0;\C)=\C$ and $H^j(\U_0;\C)=0$ for $0<j<n$. When
$V$ is a $\Q$-manifold, one also has $H^n(\U_0;\C)=0$. The computation of  $H^n(\U_0;\C) \approx H_0^{n+1}(V)$
is quite difficult in general, as it may depend on the position of singularities, see
\cite{Di}, \cite{Li2}. Here $H^*_0(V)$ denotes the primitive cohomology of $V$, i.e. the cokernel of the natural monomorphism $H^*(\mathbb{CP}^{n+1}   ) \to H^*(V)$ induced by the inclusion of $V$ into $\mathbb{CP}^{n+1}$.

For $\beta \ne 1$, one can use the isomorphism $ H^m(\U_0,\LL)= H^m(F,\C)_{\beta}$, the ${\beta}$-eigenspace of the
monodromy acting on the Milnor fiber $F$ associated to $V$. In particular, $ H^m(\U_0,\LL)=0$ for $m<n$.
It is possible to construct examples such that for  $m \in \{n, n+1\}$ one has
$$ \dim  H^m(\U_0,\LL) >  \dim  H^m(\U_0,\C).$$
Indeed, consider the polynomials in \cite{Di}, p. 148, which have
a monodromy operator without the eigenvalue 1 on all the reduced
cohomology groups  $ {\tilde H}^m(F,\C)$ (equivalently, $V$ has
the same rational cohomology as $\C\PP^n$). It is not possible that
$ {\tilde H}^m(F,\C)=0$ for all $m \in \N$, by A'Campo's result on
the Lefschetz number of the monodromy, see \cite{Di2}, p. 174.
Hence there is some integer $m$ and some  $d$-root of unity $\beta
\ne 1$ such that $\dim  H^m(\U_0,\LL_{\beta}) > 0= \dim
H^m(\U_0,\C)$. Using the Euler characteristic equality  $\chi(\U_0,\C)=\chi(\U_0,\LL_{\beta}    )$, it
follows that the inequality should hold for the two possible
values of $m$. By the minimality property of hyperplane
arrangement complements, it is known that the above inequality is
impossible for such complements, \cite{DP}.

\end{example}
The discussion in the previous section  relating local systems to connections can be extended to this setting in an obvious way. For instance, we should use now the 1-form
$$ \omega _{ \alpha  }=\sum _{j=1,s} \alpha _j \frac{df_j}{f_j}$$
where the residues $ \alpha $ satisfy the condition  $ \sum _{j=1,s} d_j \cdot \alpha _j=0$, which is a necessary condition in order to have a 1-form on $\U_0$.

\subsection{Some 2-component Arrangements}

We consider now in detail the case of hypersurface arrangements $V$ with $s=2$ irreducible components. We assume moreover that:

\noindent (i) $n>1$ and each $V_i$ has at most isolated singularities and is a $\Q$-manifold;

\noindent (ii) $V'=V_1 \cap V_2$  has at most isolated singularities; this condition is automatically fulfilled when $d_1<d_2$ and $V_2$ is smooth, see  \cite{CD}.

Let $\U_i=\mathbb{CP}^{n+1} \setminus V_i$. Then the Mayer-Vietoris sequence of the covering
$\U'=\U_1 \cup \U_2$ reads like
\begin{equation} \label{MV}
... \to H^{k-1}(\U_0) \to H^k( \U'   ) \to H^k(\U_1) \oplus H^k(\U_2)\to H^{k}(\U_0)\to ...
\end{equation}
Here and in the sequel the constant coefficients $\C$ are used unless stated otherwise.
Using Example \ref{onecomp} to handle the cohomology groups  $H^*(\U_i)$ for $i=1,2$ and the Alexander duality
isomorphism (which is compatible with the MHS due to the Tate twist $(-n-1)$)
\begin{equation} \label{AD}
H^{k}(\U')=H^{2n+2-k}(\mathbb{CP}^{n+1},V')^{\vee}(-n-1) =H^{2n+1-k}_0(V')^{\vee}(-n-1)
\end{equation}
we get the following result.

\begin{prop} \label{codim2} With the above notation and assumptions, the following hold.

\noindent (i) $H^0(\U_0)=\C$ is pure of type (0,0) and  $H^1(\U_0)=\C$ is pure of type (1,1) and is spanned by
the 1-form
$$ \omega _1=d_2 \cdot  \frac{df_1}{f_1}-d_1  \cdot  \frac{df_2}{f_2}.  $$

\noindent (ii) $H^k(\U_0)=0$ for $1<k<n$.

\noindent (iii)  $H^n(\U_0)$ is pure of weight $n+2$ and $b_n(\U_0) \leq \dim H^n_0(V')$. Moreover $H^n(\U_0)=0$ if $d_1<d_2$ and $V_2$ is smooth.

\noindent (iv)  $H^{n+1}(\U_0)$ has weights $\geq n+2$ and one has an isomorphism of MHS
$$H^{n+1}(\U_0)/W_{n+2}H^{n+1}(\U_0)= H^{n-1}_0(V')^{\vee}(-n-1).$$

\end{prop}

\proof

The vanishing of $H^n(\U_0)$ in the third claim follows from an
unexpected source. Indeed, the Gysin sequence of the smooth
divisor $X_2=V_ 2\cap \U_1$ in $\U_1$ gives a
monomorphism $H^n(\U_0) \to H^{n-1}(X_2) $. But this latter group $ H^{n-1}(X_2) $ is
trivial by some general connectivity results recently obtained by the first author, see \cite{Dtop}. The examples given in  \cite{Dtop} show that the case $d_1=d_2$ is
much more complicated, in particular the group $ H^{n-1}(X_2) $ can be non-zero.  The example below shows that the assumption $V_2$
smooth cannot be relaxed to $V_2$ with isolated singularities and
a $\Q$-manifold. The key point there is that the singularities of $V_2$ are situated on $V_1$,
a situation not covered by the results in \cite{Dtop}.

The only other claims that are not obvious are those on the MHS.
They follow from the fact that $H^n_0(V')$ has a pure HS of weight
$n$ (the singularities of $V'$ being isolated) and the following
consequence of the Alexander duality (\ref{AD})
\begin{equation} \label{Hn}
h^{p,q}(H^{k}(\U'))=h^{n+1-p,n+1-q}(H^{2n+1-k}_0(V')).
\end{equation}

\endproof

For a rank one local system $\LL \in \LL(\U_0)$, we can choose the corresponding form
$ \omega _{ \alpha  }$ to be a multiple $a(\alpha) \omega _1$ of the 1-form $ \omega _1$
introduced above. Then Propositions \ref{proppure} and \ref{codim2} yield the following.

\begin{cor} \label{codim2cor} For a non-trivial rank one local system $\LL \in \LL(\U_0)$
for which an admissible choice of residues $\alpha=( d_2 \cdot a(\alpha) , -d_1 \cdot a(\alpha))$
exists,
the following hold.

\noindent (i) $H^k(\U_0,\LL)=0$ for $k<n$.

\noindent (ii) If $ H^n_0(V')=0$ or if $d_1<d_2$ and $V_2$ is smooth, then  $H^n(\U_0,\LL)=0$.

\end{cor}

Note that the first claim above holds by Corollary \ref{cor1}, since $V'$ has only INNC singularities.

The vanishing of  $ H^n_0(V')$ holds when $V'$ is a $\Q$-homology manifold, but also in many other cases,
see for instance the discussion in \cite{Di}, pp. 207-216. There one considers only the case when
$V_1$ is a hyperplane. Indeed, any hypersurface $W$ having only isolated singularities in $\mathbb{CP}^{n}$
can be obtained as the intersection of a smooth hypersurface $V_2$ in $\mathbb{CP}^{n+1}$
with the hyperplane $H=\mathbb{CP}^{n}$, see \cite{Dif}, p. 206. However, this situation is usually uninteresting according to the second claim of the above corollary.

We conclude with an example where $V_1$ is a hyperplane, and $V_2$
is singular, so that it may have been considered already in the
previous section (in such a case $\U_0$ from this section is
exactly $\U$ from the previous section, but for the hypersurface
$V=V_2$ !).

\begin{example} \label{ex3} \rm

In  $\mathbb{CP}^3$ (with homogeneous coordinates $(x:y:z:t)$)
consider the hyperplane $H=V_1:t=0$ and the surface $V_2:
xyz-t^3=0$. Then $V_2$ has exactly 3 singularities of type $A_2$,
hence it is a $\Q$-manifold. Moreover, $H$ is transversal to $V_2$,
except at the 3 singular points of $V_2$.

\medskip

To compute the cohomology of the
complement $\U_0$, we use the Gysin exact sequence of the smooth
divisor $D=V_2 \setminus V_1$ in the affine space $\mathbb{CP}^3
\setminus V_1$ (with coordinates $(x,y,z)$) and get
$$H^k(\U_0) = H^{k-1}(D)(-1)$$
for $k=2,3$, where $(-1)$ denotes the Tate twist. Now $D$ is given by the equation $xyz=1$, hence it is a 2-dimensional torus. It follows that

\noindent (i) $H^2(\U_0)=\C^2$ is pure of type (2,2);

\noindent (ii) $H^3(\U_0)=\C$ is pure of type (3,3).
Moreover, as explained in the forth section, the 1-form $\omega _1$ is a multiple of $\frac{dg}{g}$ with $g=xyz-1$.

Let $g_0=g+1=xyz$, and note that $F_{g_0}=D$ is the Milnor fiber
of the homogeneous polynomial $g_0$. Since $\U_0=\C^3 \setminus
F_{g_0}  $, we may use the description of the cohomology groups of
$\U_0$ using Remark (2.11) in \cite{Di}, p.192. By taking
$$\eta _i =A_ixdy \wedge dz -B_idx  \wedge dz +C_i dx   \wedge dy$$
in the formula (2.12) loc.cit. with $(A_1,B_1,C_1)=(x,-y,0)$ and  $(A_2,B_2,C_2)=(0,-y,z)$ we get a basis of
$H^2(\U_0)$. A direct computation then shows that $\omega _1 \wedge \eta _i=0$ in $H^3(\U_0)=H^2(D)$.
To see this, note that $d \eta _i=dg \wedge  \eta _i=0$ and hence the Poincar\'e-Leray residue of the form
$$ \omega _1 \wedge \eta _i= \frac{dg}{g} \wedge \eta _i$$
is the form $\eta _i$. Since $d \eta _i=0$ on $\C^3$, it follows the $ \eta _i=
d \eta _i'$, for some 1-forms $\eta _i'$ on $\C^3$. Hence the cohomology class of
$\eta _i$ in $H^2(D)$ is trivial.

  It follows that, for a non-trivial rank one local system $\LL \in \LL(\U_0)$
for which an admissible choice of residues $\alpha$
exists, one has $H^*(\U_0, \LL)=H^*(H^*(\U_0), \omega _{\alpha}).$ Therefore we get the following equalities.
$$\dim H^0(\U_0, \LL) = \dim H^1(\U_0, \LL) =0, ~~\dim H^2(\U_0, \LL) =2, ~~
 \text{and} \dim H^3(\U_0, \LL)=1.$$
 In particular, $\text{Supp}  A^2(\U_0)$
coincides to the character torus $\T^1$. Indeed,  $\text{Supp}  A^2(\U_0)$ is a Zariski closed subset,
with a non-empty interior by Remark \ref{admis1rk},
in the irreducible algebraic variety  $\T^1$.  The reader should compare this fact to Theorem \ref{ucam} above.

\end{example}

\providecommand{\bysame}{\leavevmode\hbox
to3em{\hrulefill}\thinspace}


\begin{thebibliography}{10}

\bibitem{A} Arapura, D., \emph{Geometry of cohomology support loci for local systems. I}.  J. Algebraic Geom.6(1997), 563--597.


\bibitem{Bo} Borel, A. et. al., \emph{Intersection cohomology}, Progress in Mathematics, vol. \textbf{50}, Birkhauser, Boston, 1984

\bibitem{Ci} Cimasoni, D., \emph{Studying the multivariable Alexander pplynomial by means of Seifert surfaces},
arXiv: math.GT/0406150

\bibitem{COG} Cogolludo-Agustin, J.I., \emph{Topological invariants of the complement to arrangements of rational plane curves}.  Mem. Amer. Math. Soc.  159  (2002),  no. 756, xiv+75 pp.


\bibitem{CO} Cohen, D.C., Orlik, P., \emph{ Arrangements and local systems}, Math. Research Letters 7(2000),299-316.


\bibitem{CS} Cohen, D. C., Suciu, A. I., \emph{Alexander invariants of complex hyperplane arrangements},
Transactions of the AMS, Vol. 351 (10), 4043-4067, (1999).


\bibitem{CS1} Cohen, D. C., Suciu, A. I.,\emph{ Characteristic varieties of arrangements}, Math.Proc.Cambridge Philos. Soc. 127(1999), 33-54.

\bibitem{CDS} Cohen, D. C,  Denham G. and  Suciu, A.I., \emph{Torsion in Milnor fiber homology,}
Algebraic and Geometric Topology 3 (2003), 511--535.

\bibitem{COD} Cohen, D.C., Dimca, A., Orlik, P., \emph{Nonresonance conditions
for arrangements}, Ann. Inst. Fourier, Grenoble, 53, 6 (2003), 1883-1896.

\bibitem{CD} Choudary, A.D.R. and Dimca, A., \emph{Hypersurface singularities, codimension two complete intersections and tangency sets,} Geom. Dedicata 24(1987), 255-260.


\bibitem{DK} Davis, J.F., Kirk, P., \emph{Lecture Notes in Algebraic
Topology}, Graduate Studies in Mathematics, Vol.35, AMS 2001.


\bibitem{De1} Deligne, P., \emph{Equations diff\'erentielles \`a
    points singuliers r\'eguliers,} Lecture Notes in
    Math., \textbf{163}, Springer, Berlin (1970).


\bibitem{De2}  Deligne, P., \emph{Th\'eorie de Hodge II.} Publ. Math. IHES,\textbf{ 40}, 5--57 (1972).


\bibitem{Dif} Dimca, A., \emph{Topics on Real and Complex Singularities,} Vieweg Advanced Lecture in
Mathematics, Friedr. Vieweg und Sohn, Braunschweig, 1987.

\bibitem{Dtop}  Dimca, A.,  \emph{On the connectivity of some complete intersections,} arXiv:math.AG/0507501.


\bibitem{Di} Dimca, A., \emph{Singularities and Topology of Hypersurfaces}, Universitext, Springer-Verlag, 1992

\bibitem{Di1} Dimca, A.,\emph{ Hyperplane arrangements, M-tame polynomials and twisted cohomology}, in:
Commutative Algebra, Singularities and Computer Algebra,  Eds. J. Herzog, V. Vuletescu, NATO Science Series, Vol. 115, Kluwer 2003, pp. 113-126.


\bibitem{Di2} Dimca, A., \emph{Sheaves in Topology}, Universitext, Springer-Verlag, 2004.

\bibitem{DiL} Dimca, A., Lehrer, G.I., \emph{ Purity and equivariant weight polynomials,} in: Algebraic Groups and Lie Groups, editor G.I.Lehrer, Cambridge University Press, 1997.

\bibitem{Di3} Dimca, A., Nemethi, A., \emph{Hypersurface complements, Alexander
modules and Monodromy}, Proceedings of the 7th Workshop on Real and Complex Singularities,
Sao Carlos, 2002, Amer.Math.Soc (2004).

\bibitem{Di4} Dimca, A., Libgober, A., \emph{Regular functions transversal at infinity},
arXiv: math.AG/0504128

\bibitem{DiLi} Dimca, A., Libgober, A., \emph{Local topology of reducible
divisors}, arXiv: math.AG/0303215

\bibitem{DP} Dimca, A., Papadima, S., \emph{Hypersurface complements, Milnor fibers and higher homotopy groups
of arrangements},
Annals of Math. 158 (2003), 473--507.


\bibitem{EN} Eisenbud, D., Newmann, W., \emph{Three-dimensional link theory and invariants of plane curve
singularities}, Annals of Math. Studies 110, Princeton University
Press, Princeton, 1985.

\bibitem{ESV} Esnault, H., Schechtman, V., Viehweg, E.: \emph{
Cohomology of local systems on the complement of hyperplanes.}
Invent. Math.,  \textbf{109}, 557--561 (1992).
Erratum, ibid. \textbf{112}, 447 (1993)

\bibitem{GM} Goresky, M., MacPherson, R., \emph{Intersection homology II}, Invent. Math. \textbf{72} (1983), 77-129


\bibitem{Ha} Hartshorne R., \emph{Algebraic Geometry}, GTM
\textbf{52}, Springer 1977.

\bibitem{Hat} Hatcher, A., \emph{Algebraic Topology}, Cambridge
University Press, 2002.

\bibitem{Hi} Hillman, J. A., \emph{Alexander ideals of links}, LNM
\textbf{895}, Springer 1981.

\bibitem{L} Levine, J., \emph{Knot Modules, I}, Transactions of
the A.M.S., 229(1977), 1-50.


\bibitem{Li2} Libgober, A., \emph{Homotopy groups of the complements to singular hypersurfaces, II},
Annals of Mathematics, \textbf{139} (1994), 117-144

\bibitem{Li3} Libgober, A., \emph{Hodge decomposition of Alexander invariants} Manuscripta Math., \textbf{107} (2002), 251-269

\bibitem{Li6} Libgober, A., \emph{On the homology of finite abelian
covers}, Topology and its applications, \textbf{43} (1992)
157-166.

\bibitem{Li7} Libgober, A., \emph{Characteristic varieties of algebraic curves}
arXiv: math.AG/9801070, in: C.Ciliberto et al.(eds), Applications of Algebraic
Geometry to Coding Theory, Physics and Computation, 215-254, Kluwer, 2001.

\bibitem{Li8} Libgober, A., \emph{Isolated non-normal crossing},
in Real and Complex Singularities, 145-160, Contemporary
Mathematics, \textbf{354}, 2004

\bibitem{Li9} Libgober, A., \emph{Homotopy groups of complements to ample
divisors}, arXiv: math.AG/0404341, to appear in Proceedings of 12th MSJ-IRI symposium
"Singularity theory and its applications" Hokkaido Univerisity.



\bibitem{LY} Libgober, A., Yuzvinski, S., \emph{Cohomology of the Orlik-Solomon algebras and local systems}, Compositio Math. 121 (2000),337-361.


\bibitem{LV} Loeser, F and Vaqui\'e, M.:  \emph{Le polyn\^ome d'Alexander d'une courbe plane projective,} Topology 29(1990),163-173.

\bibitem{Ma} Massey, D. B., \emph{Introduction to perverse sheaves and vanishing cycles} in \emph{Singularity Theory}, ICTP 1991, Ed. D.T.Le, K.Saito, B. Teissier, 487-509


\bibitem{MS}  Matei, D. and  Suciu, A.,
{\em Hall invariants, homology of subgroups, and characteristic varieties},
Internat. Math. Res. Notices \textbf{9} (2002), 465--503.



\bibitem{Max} Maxim, L., \emph{Intersection homology and Alexander modules of
hypersurface complements}, arXiv: math.AT/0409412, to appear in Comm. Math. Helv.

\bibitem{Mi} Milnor, J., \emph{Singular points of complex hypersurfaces}, Annals of Mathematical Studies 61, vol. 50, Princeton University Press, Princeton, 1968.

\bibitem{Mi2} Milnor, J., \emph{Infinite cyclic coverings}, Topology of
Manifolds, Boston 1967.

\bibitem{Mi3} Milnor, J., \emph{A duality theorem for Reidemeister
torsion}, Ann. of Math. 76 (1962), 137-147.

\bibitem{Ra} Randell, R.:  \emph{Milnor fibers and Alexander polynomials of plane curves,} Proc. Symp. Pure Math.
40, Part 2., AMS 1983, pp. 415-419.


\bibitem{STV} Schechtman, V., Terao, H., Varchenko, A.:
 \emph{Local systems over complements of hyperplanes and the
Kac-Kazhdan condition for singular vectors.}
J. Pure Appl. Algebra, \textbf{100}, 93--102 (1995)



\bibitem{Sh} Schurmann, J., \emph{Topology of Singular Spaces and Constructible
Sheaves}, Birkhauser, Monografie Matematyczne \textbf{63}, 2003


\bibitem{Su} Suciu, A. I., \emph{Translated tori in the characteristic varieties of complex hyperplane arrangements}. Arrangements in Boston: a Conference on Hyperplane Arrangements (1999).  Topology Appl.  118  (2002),  no. 1-2, 209--223.


\bibitem{S} Sumners, D.W., Woods, J.M., \emph{The monodromy of reducible plane
curves}, Inventiones Math. \textbf{40}, 107-141 (1977).


\bibitem{W} Weibel, C.A.: \emph{An Introduction to Homological Algebra}, Cambridge Studies in Advanced Math., \textbf{38}, Cambridge Univ. Press, Cambridge (1994).

\end{thebibliography}
\end{document}